\documentclass[%
aip,
amsmath,amssymb,
preprint,%
]{revtex4-1}

\usepackage[english]{babel}

\usepackage{graphicx}
\usepackage{dcolumn}
\usepackage{bm}

\usepackage{multirow}
\usepackage{makecell}
\setcellgapes{4pt}

\usepackage[utf8]{inputenc}
\usepackage[T1]{fontenc}
\usepackage{mathptmx}

\usepackage[usenames]{color}

\usepackage[colorlinks=true]{hyperref}
\hypersetup{urlcolor=blue, citecolor=blue, linkcolor=blue}

\makeatletter
\DeclareRobustCommand*\textsubscript[1]{%
  \@textsubscript{\selectfont#1}}
\def\@textsubscript#1{%
  {\m@th\ensuremath{_{\mbox{\fontsize\sf@size\z@#1}}}}}
\makeatother

\begin{document}





\title{The limits of sustained self-excitation and stable periodic pulse trains in the Yamada model with delayed optical feedback}


\author{Stefan Ruschel}
\email{stefan.ruschel@auckland.ac.nz}
\affiliation{%
Department  of  Mathematics,  The  University  of  Auckland,  New Zealand
}%
\affiliation{%
Dodd-Walls  Centre  for  Photonic  and  Quantum  Technologies, New Zealand
}%

\author{Bernd Krauskopf}
\email{b.krauskopf@auckland.ac.nz}
\affiliation{%
Department  of  Mathematics,  The  University  of  Auckland,  New Zealand
}%
\affiliation{%
Dodd-Walls  Centre  for  Photonic  and  Quantum  Technologies, New Zealand
}%

\author{Neil G. R. Broderick}
\email{n.broderick@auckland.ac.nz}
\affiliation{%
Dodd-Walls  Centre  for  Photonic  and  Quantum  Technologies, New Zealand
}%
\affiliation{%
Department  of  Physics,  The  University  of  Auckland,  New Zealand
}%

\date{\today}

\begin{abstract}
We consider the Yamada model for an excitable or self-pulsating laser with saturable absorber, and study the effects of delayed optical self-feedback in the excitable case. More specifically, we are concerned with the generation of stable periodic pulse trains via repeated self-excitation after passage through the delayed feedback loop, as well as their bifurcations.
We show that onset and termination of such pulse trains correspond to the simultaneous bifurcation of countably many fold periodic orbits with infinite period in this delay differential equation. We employ numerical continuation and the concept of reappearance of periodic solutions to show that these bifurcations coincide with codimension-two points along families of connecting orbits and fold periodic orbits in a related advanced differential equation. These points include heteroclinic connections between steady states, as well as homoclinic bifurcations with non-hyperbolic equilibria. Tracking these codimension-two points in parameter space reveals the critical parameter values  for the existence of periodic pulse trains. We use the recently developed theory of temporal dissipative solitons to infer necessary conditions for the stability of such pulse trains. 
\end{abstract}

\maketitle

\begin{quotation}
	The generation of periodic light pulses constitutes the basis for modern telecommunication\cite{Agrawal2002}, material processing\cite{Steen2010} and high energy physics\cite{Gibson1974}. One particular mechanism that allows for the formation of high-intensity pulses is Q-switching in semiconductor lasers with saturable absorber \cite{Ueno1985,Erneux1988,Yamada1993}. Here, the absorber section acts as a situational shutter that prevents low-intensity light from leaving the laser.
	The laser is effectively off, yet a  sufficiently large perturbation is able to cause a brief, high-intensity pulse of light\cite{Dubbeldam1999,Dubbeldam1999a}. Noise present in the laser may lead to sporadic pulsation, and in this work we are concerned with a pace-making strategy, where the output pulse is re-injected into the laser cavity to trigger a subsequent pulse\cite{Otto2012,Jaurigue2015}. In order to understand this feedback mechanism, we consider the Yamada model\cite{Yamada1993} subject to delayed optical self-feedback, and study periodic solutions that are close to homoclinic in a corresponding (advanced) delay differential equation \cite{Lin1986, Lin1990}. In this way, we are able to characterize the model parameters that allow for stable self-pulsation generated by the effect of delayed feedback on the excitable laser.
	
	Our methodology  is not limited to lasers, but is relevant for the rigorous analysis of delay-coupled excitable systems in general. Future applications will include recurrent networks of pulsing lasers or neurons, which are in a similar excitable configuration when at rest\cite{Izhikevich2005}; as well as the exploration of this analogy\cite{Garbin2017} with respect to neuro-morphic computation tasks performed by coupled laser systems\cite{Bueno2018,Selmi2015,Prucnal2016}.\\
\end{quotation}


\section{Introduction}

Lasers display a wide range of nonlinear dynamical phenomena beyond the steady output of light, ranging from excitable, to oscillatory and chaotic behavior \cite{Krauskopf2000, Ludge2011}. A well known way to achieve this dynamical richness are time delays caused by optical feedback\cite{RuschelYanchuk2017, Soriano2013, YanchukGiacomelli2017}. Their high-dimensional effects are exploited in many applications, such as private communication\cite{Argyris2005}, cryptography\cite{Ludge2011} and photonic computing \cite{Appeltant2011, Pammi2019}. 

This work is motivated by experiments with a micro-pillar semiconductor laser with a saturable absorber and a reflecting mirror\cite{Terrien2017a,Terrien2018a,Terrien2019}. We study how the interplay of the absorbing medium and the delayed feedback allows for periodic pulse trains, where one or more temporally localized pulses of light are sustained in the optical feedback loop. We follow up on earlier work\cite{Walker2008,Terrien2017,Terrien2018a,Terrien2019}, and study this set up theoretically via the delay differential equation 
\begin{align}
G^{\prime}(t)= & \gamma_G\left(A-G(t)[1+I(t)]\right),\label{eq:G-def}\\
Q^{\prime}(t)= & \gamma_Q\left(B-Q(t)[1+aI(t)]\right),\label{eq:Q-def}\\
I^{\prime}(t)= & \left[G(t)-Q(t)-1\right]I(t)+\kappa I(t-\tau).\label{eq:I-def}
\end{align}
It describes the time evolution of the laser output intensity $I(t)\geq0$, loosely speaking the amount of photons emitted by the laser, as well as the gain $G(t)\geq0$ and absorption $Q(t)\geq0$ of photons in the cavity of the micro-pillar laser\cite{Terrien2018}. Equations~(\ref{eq:G-def})--(\ref{eq:I-def}) reflect the influence of delayed optical feedback on a solitary laser ($\kappa=0$) described by the Yamada model\cite{Yamada1993} in dimensionless form\cite{Dubbeldam1999}, which is a prototypical model for the creation of high-intensity pulses in lasers with saturable absorber via a specific mechanism called Q-switching \cite{Ueno1985,Erneux1988,Yamada1993, Barbay2011}. 

We proceed with the explanation of the model. From Eq.~(\ref{eq:I-def}), one immediatly sees that $I(t)$ increases fast when $G(t)>Q(t)+1$, which corresponds to the amplification of light in the laser cavity.  This threshold behavior and the maximum value of the gain $G(t)$ are largely determined by the so-called pump parameter $A$, which is the main control parameter in the experiment, as well as the specific properties of the absorbing medium described by the parameters $B,a,\gamma_G,$ and $\gamma_Q$. Here, $B$ and $a$ characterize the loss of photons in the absorber section; $\gamma_G$ and $\gamma_Q$ correspond to the attenuation rates of photons in the gain and absorber sections; and $0\leq\kappa\leq1$ and $\tau\geq0$ are strength and delay of the optical feedback, respectively. We restrict ourselves to the specific choice of parameters $B=1.8,$ $a=1.8,$ and $\gamma_G=\gamma_Q=\gamma=0.04$ that has been considered before\cite{Dubbeldam1999,Walker2008,Terrien2017}, and study behavior of system~(\ref{eq:G-def})--(\ref{eq:I-def}) with respect to the parameters $A,\kappa$ and $\tau$. 

We emphasize that the phase of the laser light does not play a role here, as the feedback amplitude is the only determining factor for whether a subsequent pulse is triggered or not. The general intuition is as follows: in between pulses the laser is off, and there is no well defined phase in the cavity, as there is no amplification of light. In the same way, the response of the laser does not depend on the phase of the perturbation.

\begin{figure}[]
	\centering{}\includegraphics[width=.6\linewidth]{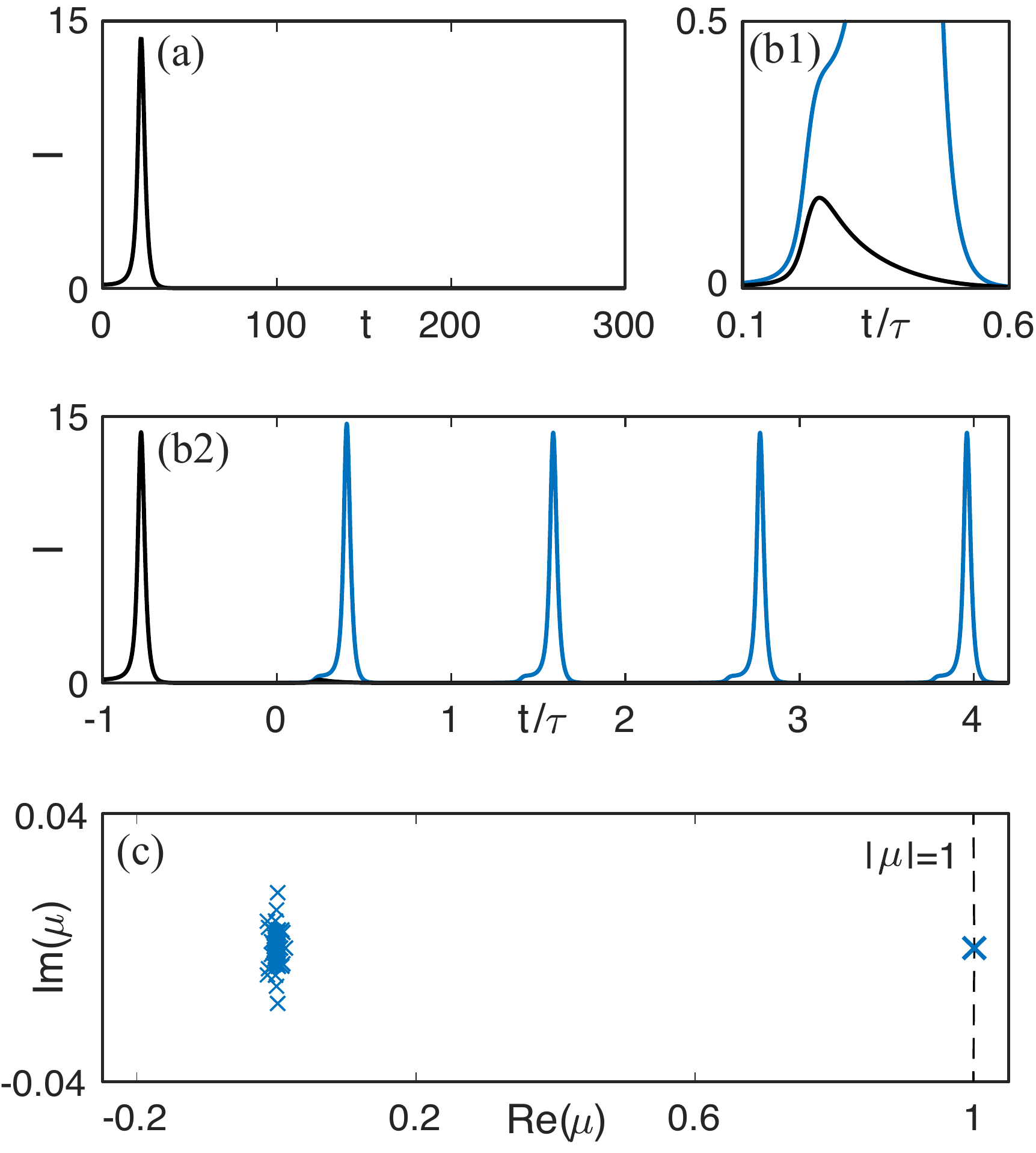}
	\caption{\label{fig:1}
		(Color online) Pulse genesis and creation of a stable periodic pulse train. (a) Single large pulse of Eqs.~(\ref{eq:G-def})--(\ref{eq:I-def}) without feedback $(\kappa=0)$; (b1) response of Eqs.~(\ref{eq:G-def})--(\ref{eq:I-def}) to initial pulse in (a) for $\kappa=0.005$ (black) and $(\kappa=0.01)$ (blue);  (b2) periodic pulse train (blue) for $(\kappa=0.01)$; (c) Floquet spectrum of the periodic pulse train in (b2). Other parameters are $A=6.5$, $B=5.8$, $a=1.8$, $\gamma_G=\gamma_Q=0.04$, $\tau=100$. 
	}
\end{figure}

To start off, we briefly review the bifurcation structure of Eqs.~(\ref{eq:G-def})--(\ref{eq:I-def}) without feedback, i.e. for $\kappa=0$ when the system is an ODE; the details can be found in Refs.~\onlinecite{Yamada1993,Dubbeldam1999,Huber2005, Otupiri2019}. For the considered parameter values and small values of $A$, the laser is off, and the corresponding steady state \textsf{o} (the "off state"), given by $(G,Q,I)=(A,B,0)$, is globally exponentially stable. Increasing $A$, we encounter a fold (or saddle-node) bifurcation that gives rise to two unstable steady states \textsf{p} and \textsf{q}, in addition to the still stable off state \textsf{o}. The steady states \textsf{p} and \textsf{q} are of the form $(G,Q,I)=(A/(1+I),B/(1+aI),I),$ where $I_\mathsf{p}<I_\mathsf{q}$ are the two solutions of the algebraic equation $0=A/(1+I)-B/(1+aI)-1.$ More specifically, \textsf{p} is a saddle with a one-dimensional unstable and a two-dimensional stable manifold that locally separates the phase space near \textsf{p}. This is an excitable configuration:
the laser is essentially off, yet a sufficiently large perturbation that brings the system across the stable manifold of \textsf{p}, causes a fast increase in $I$ that is followed by fast relaxation as the solution follows the unstable manifold of \textsf{p} back to \textsf{o}; see Refs.~\onlinecite{Dubbeldam1999, Otupiri2019} for details. Such a trajectory corresponds to the emission of a single pulse of light by the laser; see Fig.~\ref{fig:1}(a). 

One can imagine that resetting the value of $I$ at some point $t$, we can trigger another pulse. In a solitary laser such perturbations arise due to noise fluctuations in the laser cavity, and result in sporadic pulsation of the laser \cite{Dubbeldam1999a, Terrien2017a}. A similar role can be played by the delayed feedback. 
Consider the single pulse of the system without feedback as the initial condition of the full system (\ref{eq:G-def})--(\ref{eq:I-def}). 
At time $t=0$, the pulse has 'traveled' through the delay loop and acts as a perturbation to the laser off-state; if it is sufficiently large, it should trigger a subsequent pulse. This is illustrated in Fig.~\ref{fig:1}, where panel (a) shows the input pulse. 
Figure~\ref{fig:1}(b1) shows the response  of system~(\ref{eq:G-def})--(\ref{eq:I-def}) to the input pulse for the two different values of the coupling strength $\kappa=0.005$ (black) and $\kappa=0.01$ (blue). Whereas the response is small for $\kappa=0.05$ (black), the choice $\kappa=0.01$ (blue) leads to a large secondary pulse.  
This secondary pulse triggers yet another pulse after having traveled the delay loop, and so forth, thus establishing a stable periodic pulse train generated by repeated self-excitation; see Fig.~\ref{fig:1}(b2).
Figure~\ref{fig:1}(c) shows the numerically computed Floquet spectrum of the pulse train shown in Fig.~\ref{fig:1}(b2) with all Floquet multipliers well within the unit circle; see for example Refs.~\onlinecite{Yanchuk2009, Sieber2013a, Yanchuk2019} for details on the Floquet spectrum of periodic solutions in delay differential equations. 

It is the main goal of this article to determine, independently of the initial condition, the minimum coupling strength $\kappa_{\min}=\kappa_{\min}(A,B,a,\gamma_G, \gamma_Q)$, that gives rise to stable periodic pulse trains, as well as to determine and the corresponding bifurcations in Eqs.~(\ref{eq:G-def})--(\ref{eq:I-def}).  
In addition, we discuss the bifurcations involved in the termination of periodic pulse trains due to either strong contraction towards the constant lasing sate \textsf{q} (after Hopf bifurcation through a corresponding change of parameters), or loss of stability of the off-state \textsf{o}. In particular, we also compute the corresponding maximum coupling strength $\kappa_{\max}=\kappa_{\max}(A,B,a,\gamma_G, \gamma_Q)$. 
We obtain these results by numerical continuation of periodic solutions and their bifurcations in Eqs.~(\ref{eq:G-def})--(\ref{eq:I-def}) with the continuation package DDE-Biftool \cite{Engelborghs2002, Sieber2014}. 

The article is structured as follows. Section~\ref{sec:results} outlines our main results; it contains a numerically computed diagram of the critical parameter values corresponding to the onset and termination of periodic pulse trains. In Sec.~\ref{sec:methods}, we introduce the concept of reappearance of periodic solutions \cite{Yanchuk2009}, which allows us to identify families of periodic pulse trains with families of periodic solutions close to a homoclinic bifurcation in a related (advanced) delay differential equation. Section~\ref{sec:methods} also contains a brief introduction to homoclinic orbits in (advanced) delay differential equations\cite{Lin1986, Hupkes2009}, and the related concept of temporal dissipative solitons\cite{Yanchuk2019}, which is used to obtain necessary conditions for the stability of a pulse train for sufficiently large values of the delay. Section~\ref{sec:2parametercont} discusses two-parameter bifurcation diagrams. We encounter codimension-two points along families of homolinic orbits and fold periodic orbits that are organizing centers for the existence of periodic pulse trains; projecting these points onto the $(A,\kappa)$-parameter plane reveals the critical coupling strengths with respect to the pump parameter $A$. Section~\ref{sec:discussion} discusses our findings with respect to the experiment. 

\section{Overview of main results}\label{sec:results}

\begin{figure}[]
	\centering{}
	\includegraphics[width=.6\linewidth]{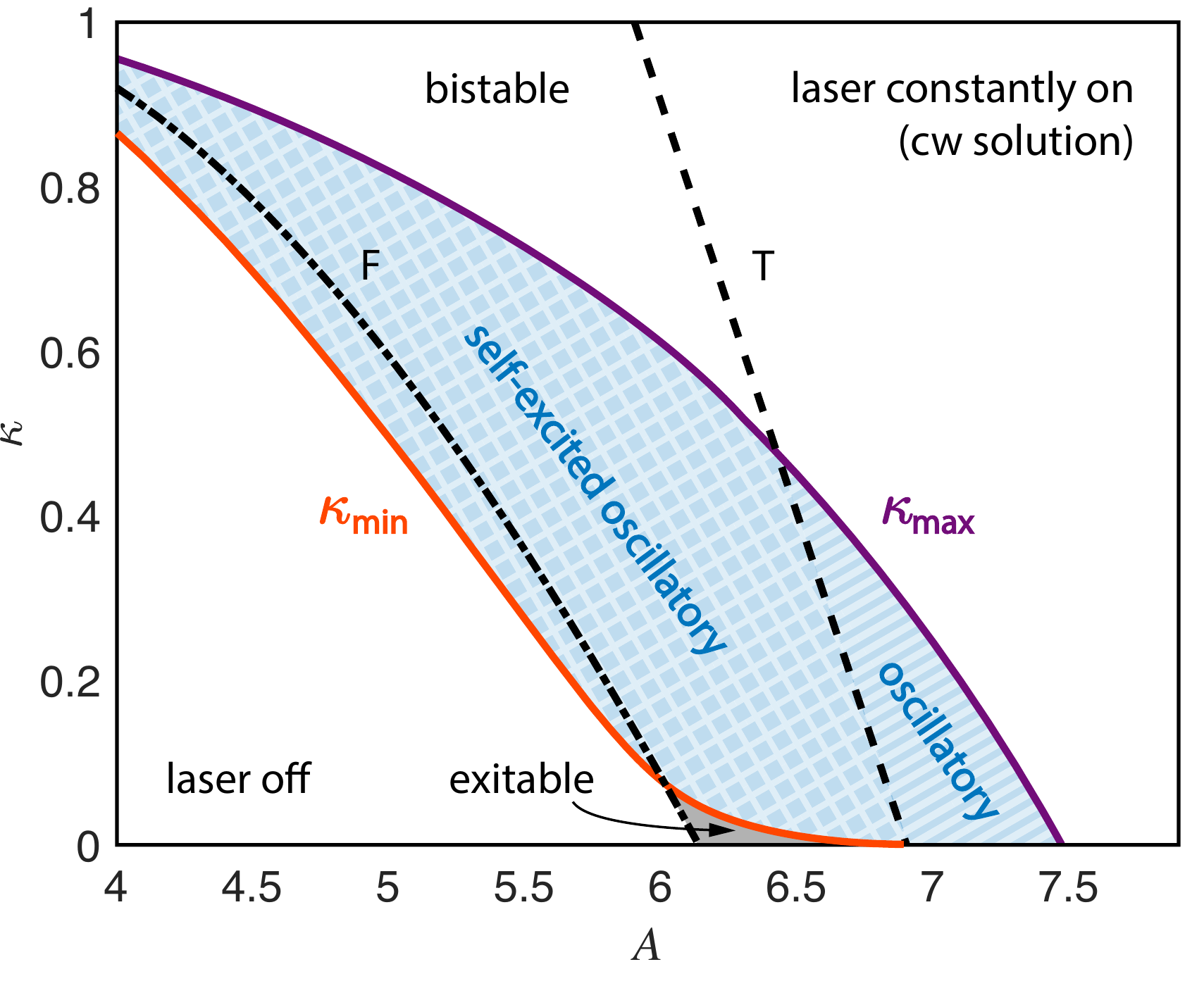}
	\vspace{-.6cm}
	\caption{(Color online) Parameter region of periodic pulse trains (blue). The panel shows a projection of the $(A,\kappa,\tau)$-bifurcation diagram onto $(A,\kappa)$-plane. The curves of minimum (orange) and maximum (purple) coupling strength correspond to codimension-two points along families of homoclinic orbits or fold perioidic orbits. Dashed curves correspond to delay independent codimension-one bifurcations; see the text for details. Other parameters are $a=1.8$, $B=5.8$, $\gamma_G=\gamma_Q=0.04$. \label{fig:2}
	}
\end{figure}

This section summarizes our main results. Figure~\ref{fig:2} shows the minimum coupling strength $\kappa_{\min}(A;B,a,\gamma_G, \gamma_Q)$ (orange) and the maximum coupling strength $\kappa_{\max}(A;B,a,\gamma_G, \gamma_Q)$ (purple) that allow for periodic pulse trains in system (\ref{eq:G-def})--(\ref{eq:I-def}) with fixed parameters $a,B,\gamma_G,\gamma_Q$. Figure~\ref{fig:2} is not obtained by direct simulation and comparison of the time series, but shows a projection of the $(A,\kappa,\tau)$-bifurcation diagram onto the $(A,\kappa)$-parameter plane; it is therefore independent of the initial condition. We obtain $\kappa_{\min}$ and $\kappa_{\max}$ by tracking codimension-two points along families of homoclinic orbits and fold periodic orbits of Eqs.~(\ref{eq:G-def})--(\ref{eq:I-def}); see Sec.~\ref{sec:2parametercont} for details. In addition, Fig.~\ref{fig:2} displays the fold bifurcation curve \textsf{F} (dashed dotted)  of the steady states \textsf{p} and \textsf{q} given by
$$
\kappa_\mathsf{F}(A;B,a)=\frac{(1-A)a-B-1+2\sqrt{aAB}}{a-1},
$$
and the transcritical bifurcation curve \textsf{T} (dashed)  of the steady states \textsf{o} and \textsf{p} given by
$$ 
\kappa_\mathsf{T}(A;B)=-(A-B-1),
$$
as parametrized by $A$; both $\kappa_\mathsf{F}$ and $\kappa_\mathsf{T}$ are independent of the feedback delay. To the left of the fold curve \textsf{F} the laser off state \textsf{o} exists, and there are no further steady states; between the curves \textsf{F} and \textsf{T}, three steady states \textsf{o},\textsf{p},\textsf{q} exist; and to the right of \textsf{T} only the laser off state \textsf{o} and the constant amplitude lasing state \textsf{q} (also called continuous wave solution) are physically relevant, as the $I$-component of \textsf{p} is negative. 

The curves $\kappa_{\min}$, \textsf{F}, and the $A$-axis enclose a region (dark gray), which is an extension of the excitable region of the Yamada model without feedback to $\kappa>0$. 
Increasing $\kappa$ above $\kappa_{\min}$, one enters the parameter region of sustained self-pulsation generated via delayed self-excitation (blue checkered). We call this regime self-excited oscillatory, as opposed to the oscillatory regime (blue striped) to the right of the curve  \textsf{T}  in Fig.~\ref{fig:2}. 
The oscillatory regime here is a generalization of the oscillatory regime of the system without feedback extended to $\kappa>0$. Both regimes (blue) allow for periodic solutions. 

Figure~\ref{fig:2} should be read as follows. For each $(A,\kappa)$ in the blue region, there is $\tau_0\geq0$, such that for $(A,\kappa,\tau_0)$ and fixed $a,B,\gamma_G, \gamma_Q$, there is a periodic solution that is either a periodic pulse train or an oscillation with small amplitude about \textsf{q}. 

The major distinction between the two regimes is that the self-excited oscillatory regime allows for stable periodic periodic pulse trains with arbitrarily high period $T$, where the period $T$  of a single pulse train is approximately given by the length of the delay $\tau$ plus an offset $\delta$ corresponding to the response of the laser, divided by the number $k$ of pulses per delay interval, i.e. $T\approx (\tau+\delta)/k$. 
Thus, for a periodic pulse train with $k$-pulses in a delay interval to exist and for it to be stable, it is necessary that the delay is sufficiently large. Section~\ref{sec:TDS} shows that
\begin{equation}\label{eq:crit-delay}
\tau \gg \frac{1}{1-\sqrt[k]{\frac{\kappa}{|A-B-1|}}}
\end{equation}
is a necessary condition for stability.
Conversely, the number of equidistant pulses per delay interval that give rise to a stable periodic pulse trains is bounded from above by
\begin{equation}\label{eq:crit-k}
k \ll \tau \ln\left|\frac{A-B-1}{\kappa}\right|,
\end{equation}
assuming that the delay is sufficiently large.

In the oscillatory regime, however, the maximum period of the pulse trains is bounded for arbitrarily large values of the delay; moreover, the corresponding periodic solutions are not necessarily stable for large values of the delay. This is very intuitive: To the left of the transcritical bifurcation curve \textsf{T}, the laser off state \textsf{o} is unstable (see Appendix \ref{sec:spec-o}), and the time spent close to \textsf{o} is not solely determined by the time between pulses anymore, but rather by the rate of expansion in the vicinity of the now saddle equilibrium \textsf{o}. An earlier detailed bifurcation analysis by some of the authors shows that, in the oscillatory regime, pulse trains can become increasingly unstable as the delay increases, for example via torus bifurcations\cite{Terrien2017}. This is due to a large delay instability of modulational type, which is not covered here; see Refs. \onlinecite{Yanchuk2009, Yanchuk2019, YanchukGiacomelli2017} for details. 

The regime of periodic pulse trains (blue) is bounded from above by the purple curve\linebreak $\kappa_{\max}(A;B,a,\gamma_G, \gamma_Q)$. For $\kappa>\kappa_{\max}$, there are no periodic solutions and the equilibrium \textsf{q} is stable. This corresponds to a constant amplitude output of the laser, when $\kappa>\kappa_\mathsf{T}$ (continuous wave solution), or a bistable configuration between \textsf{o} and \textsf{q}, when $\kappa<\kappa_\mathsf{T}$. 

Figure~\ref{fig:2} shows an interesting feature of system~(\ref{eq:G-def})--(\ref{eq:I-def}). Notice that the pump strength $A$ needed to create a pulse in system~(\ref{eq:G-def})--(\ref{eq:I-def}) (and therefore a pulse train) can be significantly reduced for high values of the coupling strength $\kappa$. Between the curves $\kappa_{\min}$ and \textsf{F}, the pump strength $A$ is such that the laser without feedback shows no activity; in particular, the steady state \textsf{q} that usually governs excitability is absent. Nevertheless, the delayed feedback is still able to produce a pulsed response. This can be thought of as a type of feedback-induced excitability.

\section{Basic Concepts and Methods}\label{sec:methods}
A detailed bifurcation analysis of (\ref{eq:G-def})-(\ref{eq:I-def}) was given earlier by some of the authors\cite{Walker2008,Terrien2017}. We briefly review the relevant results and focus on periodic pulse trains and their bifurcations. For the sake of accessibility, we keep technical notations to a minimum, and leave more involved computations to the Appendix. An introduction to the general theory of delay differential equations, existence
and uniqueness of solutions, and specific notions of solution and phase space will not be presented here, but can be found in classic textbooks \cite{Hale1993, Diekmann1995, Guo2013}. 

\begin{figure}[]
  \vspace{0cm}
  \begin{minipage}[c]{0.6\textwidth}
  	\vspace{0cm}
  	\centering{}\includegraphics[width=1.\linewidth]{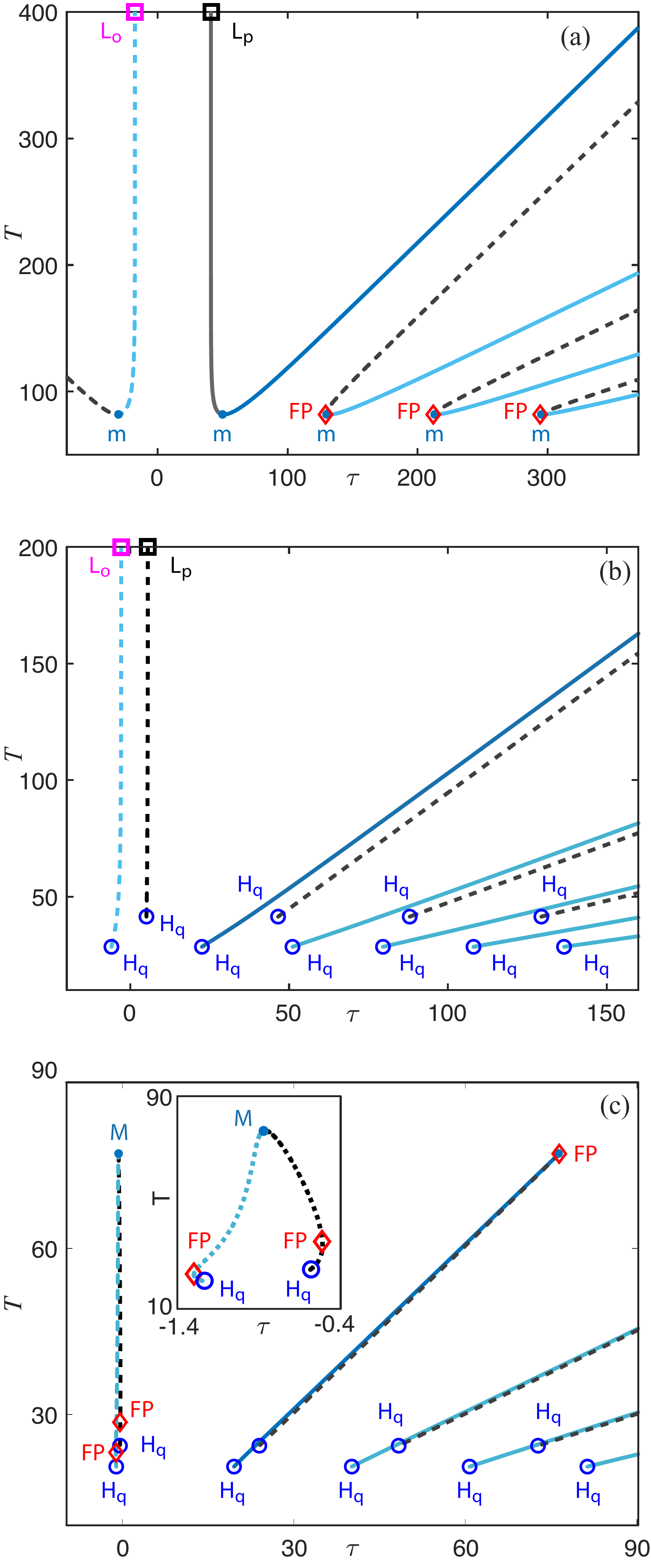}
  \end{minipage}\hfill
  \begin{minipage}[c]{0.3\textwidth}
  \vspace{-8cm}
  \caption{
  		(Color online) Bifurcation diagram of periodic solutions. Shown is the period \textsf{T}  with respect to the delay $\tau$ for different values of $\kappa$: (a) $\kappa=0.01$; (b) $\kappa=0.1$; (c) $\kappa=0.38$ with an enlargement for negative $\tau$. Other parameters are $A=6.5$, $B=5.8$, $a=1.8$, and $\gamma_G=\gamma_Q=0.04$. Shown are stable periodic solution (solid curves), unstable periodic solutions (dashed curves), homoclinic bifurcations \textsf{L\textsubscript{o}} and \textsf{L\textsubscript{p}}  (squares), Hopf bifurcations \textsf{H\textsubscript{q}}  (blue circles), fold bifurcations of periodic orbits \textsf{FP} (red diamonds) and local minima \textsf{m}  and maxima \textsf{M}  of the period with respect to the delay (small blue dots). \label{fig:3}
  	}
  \end{minipage}
\end{figure}

Following up on the introduction, we expect that the period $T$  of the pulse train in Fig.~\ref{fig:1}(b2) (and, therefore, in the self-excited oscillatory regime in Fig.~\ref{fig:2}) is slightly larger than the time, needed to complete one  delay loop. The intuition is clear: consider (\ref{eq:G-def})-(\ref{eq:I-def}) in the excitable configuration with a pulse emitted at some point in time $t$.  This pulse affects the system only after having traveled the feedback loop at time $t+\tau$. Moreover, there is a certain 'response time' $\delta$ of the system needed to produce a subsequent pulse; $\delta$ is henceforth the drift\cite{YanchukGiacomelli2017}. Consider now Fig.~\ref{fig:3}(a). Shown is the period of the pulse train in Fig.~\ref{fig:1}(b2) with $\tau=100$ (dark blue in both figures) with respect to the delay $\tau$ as a continuation parameter. Solid curves depict stable periodic solutions, and dashed curves indicate that the corresponding periodic orbit is unstable for the considered parameter values; see Sec.~\ref{sec:TDS} for details on their Floquet spectra. We observe that for large values of $\tau$, the period $T=T(\tau)$ along the branch of the periodic pulse train (dark blue) scales approximately linearly with $\tau$, and the drift $\delta(\tau)=T(\tau)-\tau>0$ approaches a constant as $\tau$ increases along the branch; we denote this constant $\delta_{0}$. Note that for large values of $\tau$, the contribution of the drift to the period $\delta(\tau)/T(\tau)\approx\delta_0/\tau$ is vanishingly small. Following the branch (dark blue) to smaller values of $\tau$, we encounter a minimum \textsf{m} of the period with respect to the delay (blue dot) such that at $\tau=\tau^\mathsf{m}_1$, we have $T(\tau^\mathsf{m}_1)=\min_\tau T(\tau)=:T_{\min}$ and $\partial_\tau T(\tau^\mathsf{m}_1)=0$, where $\partial_\tau$ is the derivative with respect to $\tau$. We color the branch segment to the left of \textsf{m} in gray. This choice is for ease of reference; it does not reflect a bifurcation, but helps with later explanations. Further decreasing $\tau$ along the solid gray branch segment, we encounter a homoclinic bifurcation point \textsf{L\textsubscript{p}}  (indicated by black square) of the steady state \textsf{p} at $\tau=\delta_1$. (Section \ref{sec:homoclinic} contains a brief introduction to homoclinic orbits in delay differential equations.) One further notices that the other branches shown in Fig.~\ref{fig:3}(a) appear to be deformed copies of the branch that we have just discussed. This property, called reappearance of periodic solutions\cite{Yanchuk2009}, is a generic feature of delay differential equations.

\subsection{Reappearance of pulse trains and local period extrema}\label{sec:reapperance}
 We show how the concept of reappearance facilitates our understanding of Fig.~\ref{fig:3}, and how Fig.~\ref{fig:3} contains additional information about periodic pulse trains for large values of $\tau$. Assume that the general delay differential equation
\begin{equation}\label{eq:DDE-general-0}
x^\prime(t)=f(x(t),x(t-\tau))
\end{equation}
with $x(t)\in\mathbb{R}^d$ for some $d\in\mathbb{N}$ and $f:\mathbb{R}^d\times \mathbb{R}^d \to \mathbb{R}^d$, possesses a periodic solution $\tilde x$ with period $T=T(\tau_0)$, that is $\tilde x(t)=\tilde x(t-T(\tau_0))$ for all $t$. Noting that $$\tilde x(t-\tau_0) = \tilde  x(t-\tau_0-kT(\tau_0)),$$ for all $k\in\mathbb{Z},$ we immediatly have that $\tilde x$ is also a solution of 
\begin{equation}\label{eq:DDE-general-k}
\tilde x^\prime(t)=f(\tilde x(t),\tilde x(t-\tau_k)),\quad\tau_k = \tau_0 + kT(\tau_0),
\end{equation}
for all $t$ and $k$. We say that the periodic solution $\tilde x$ \textit{reappears} at $\tau_k$. 
In fact, entire branches of periodic solutions parameterized by $\tau_0$ reappear via the mapping $\tau\mapsto\tau+T(\tau)$, as can be seen in Figure~\ref{fig:3}(a). Depending on the derivative of the period with respect to the delay, branch segments with $|1+kT^\prime(\tau_0)|<0,$ $k\in\mathbb{Z}$, get locally squeezed, segments with $|1+kT^\prime(\tau_0)|>0,$ $k\in\mathbb{Z}$, get locally stretched in the process\cite{Yanchuk2009}.
The solid gray and dark blue branches of periodic orbits reappear in the form of the dashed gray and cyan branches, respectively. For decreasing $\tau$ along the k-th reappearing branch for $k>1$ (cyan), we respectively encounter a fold of periodic orbits FP (red diamonds), where a real, and positive Floquet multiplier crosses the unit circle. Further destabilizing bifurcations along the gray dashed branches are possible, but are not investigated here; see Sec.~\ref{sec:TDS} for details on the stability along the blue branches. 

Consider now the stable periodic pulse train for $\tau_1=100$ of Fig.~\ref{fig:1}(b2) (on dark blue branch) and with fixed period $T(\tau_1)$ in Fig.~\ref{fig:3}(a). We call this periodic pulse train a $1$-pulse solution. The solid cyan branches in Figure~\ref{fig:3}(a) correspond to stable $k$-pulse solutions with the same period as the ones on the dark blue branch,  but accommodating $k$ pulses on the respective (longer) delay interval $\tau_k=\tau_1+(k-1)T(\tau_1)$ for $k\in\{2,3,4\}.$ In particular, for a fixed large value of the delay, the number of existing periodic pulse trains can be approximated by $j=\lfloor\tau/T_{\min}\rfloor$ (rounded down) from the $1$-pulse, to the $2$-pulse, and so on, to the $j$-pulse solution\cite{Yanchuk2009}. 
Decreasing the delay along the branch of $k$-pulse trains in Fig.~\ref{fig:3}(a), the solution  switches from $k$ pulses to $(k-1)$ pulses per delay interval close to each of the folds.  

On the other hand, the $1$-pulse solution reappears as a '$0$-pulse' solution for negative delay $\tau_0=\tau_1-T(\tau_1)<0$ (cyan dashed). Here, the delay interval is much shorter than the period of the pulse train, so that it takes multiple delay intervals to complete one period. Increasing $\tau_0$ along this branch, the period $T(\tau_0)$ increases rapidly, and we encounter a homoclinic bifurcation \textsf{L\textsubscript{o}} (indicated by a magenta square) with the \textsf{o} steady state, when $\tau_0=\tau_1-T(\tau_1)\to - \delta_0$ as $T\to\infty$. $\delta_0$ is the asymptotic drift of the $1$-pulse solution.  Note that \textsf{L\textsubscript{o}} is not directly physically relevant. It will become clear, however, that the continuation of \textsf{L\textsubscript{o}} with respect to $\kappa$ and $A$ is crucial for our understanding of the behavior of the periodic pulse train for large values of the delay. Fig.~\ref{fig:3}(a) shows that the gray branch segments reappear analogously. 

In this way, a lot of information about a $k$-pulse trains is already contained in the curve of the $0$-pulse trains. In particular, a local minimum \textsf{m}  gets mapped to a local minimum of the reappearing branch. If $\tau^{\mathsf{m}}_0=\mbox{argmin }T(\tau)$ satisfies $T(\tau^{\mathsf{m}}_0)=\min_\tau T(\tau)=:T_{\min}$, then $$T(\tau^{\mathsf{m}}_0+kT(\tau^\mathsf{m}_0))=T(\tau^\mathsf{m}_k)=T_{\min}$$ for all $k$. In addition, the period minima \textsf{m} are a good approximation of the fold points \textsf{FP}  (see Fig.~\ref{fig:3}(a)) and the same holds for fold points close to local maxima of the period with respect to the delay. 
A rigorous statement goes as follows. At a fold, the parametrization $\tau_k(\tau_0)$ is locally not invertible, leading to the condition $0=d\tau_k/d\tau_0=1+kT^\prime(\tau_0)$ for a fold at the $k$-th branch\cite{Yanchuk2009}. This  is equivalent to $T^\prime(\tau_0)=1/k$. As a result, the fold point on the $k$-th branch reappears from a point $(\tau_0,T(\tau_0))$ on the $0$-th branch with $T^\prime(\tau_0)=1/k\to0$ for $k\to0$, and there are infinitely many such bifurcations. Moreover, the fold points $(\tau_k^\mathsf{FP},T(\tau_k^\mathsf{FP}))$ converge to the period minimum $(\tau^\mathsf{m}_k,T_{\min})$ from above on the respective branch, since $$|\tau_k^\mathsf{FP}-\tau^\mathsf{m}_k|\to0$$ as $k\to\infty$. The analogous statement is true for local maxima. 

As a result, there is a direct correspondence between the homoclinic bifurcations \textsf{L\textsubscript{o}}  and \textsf{L\textsubscript{p}}  along the branches parametrized by $\tau_0$ and $\tau_1$ with the families of $k$-pulse solutions parametrized by $\tau_k$, as $\tau_k\to\infty.$ When studying $k$-pulse trains (and their fold points) we may, therefore, restrict ourselves to a finite range of delay values. 

This observation is rather powerful as it allows us to study bifurcations of periodic orbits for large values of the delay indirectly via bifurcations of the periodic and homoclinic orbits on the branch parameterized by $\tau_0$ and $\tau_1$. In particular, we can study the change of the bifurcation diagram as we decrease $\kappa$ to $\kappa_{\min}$. A two parameter continuation in the $(\tau,\kappa)$-parameter plane reveals, that as $\kappa\to\kappa_{\min}$, the minimum period $T_{\min}$ of the periodic pulse trains shown in Fig.~\ref{fig:3} grows beyond bound. At the same time, the fold periodic orbits along \textsf{FP} exist at larger and larger values of $\tau$ until they disappear simultaneously at $(\tau,\kappa)=(\infty,\kappa_{\min})$; see Sec.~\ref{sec:2parametercont} for details. 

On the other hand, we can study the qualitative change of the bifurcation diagram for larger values of $\kappa$. For $\kappa=0.1$, the two homoclinic bifurcations \textsf{L\textsubscript{o}}  and \textsf{L\textsubscript{p}}  and, therefore, the branches of $k$-pulse solutions, have moved closer together; see Fig.~\ref{fig:3}(b). Additionally, another fold of periodic orbits was created via a neutral saddle bifurcation of \textsf{L\textsubscript{o}} , and all fold periodic orbits have merged with generalized Hopf bifurcation points (Sec.~\ref{sec:2parametercont}); thereby giving rise to pairs of Hopf bifurcations \textsf{H\textsubscript{q}} , which ultimately lead to the disappearance of the local period minima and splitting of the blue (cyan) and black (gray) colored branches. For $\kappa=0.38$, the homoclinic bifurcations have disappeared entirely, which is accompanied by an appearance of infinitely many fold periodic orbits \textsf{FP} (only one is shown here) that now connect different branches of $k$-pulse trains, bounding their period from above; see Figs.~\ref{fig:3}(c). The appearance of folds marks the transition from the self-excited oscillatory regime to the oscillatory regime in Fig.~\ref{fig:2}. Further increasing $\kappa$ to $\kappa_{max}$ causes the pulse trains to dissappear entirely; see Sec.~\ref{sec:2parametercont} for details. As a remark, Fig.~\ref{fig:3}(c) is an extension of the bifurcation diagram (Fig.~6) in Ref.~\onlinecite{Terrien2017} to negative values of the delay.

\subsection{Homoclinic orbits and temporal dissipative solitons}\label{sec:homoclinic}
We now focus on the homoclinic orbits \textsf{L\textsubscript{o}}  and \textsf{L\textsubscript{p}} . In general, a homoclinic orbit of  a delay differential equation is a special, closed loop solution that converges to the same saddle equilibrium in forward and backward time. Finding such connecting orbits relies on the generalization of the corresponding stable and unstable eigenspaces of the steady state via the concept of exponential dichotomies and the existence of corresponding projection operators. This approach for delay differential equations was initiated by a seminal paper by Xiao-Bao Lin\cite{Lin1986}; and the method developed in this paper, commonly known as Lin's method\cite{Lin1990} has become the standard method for numerical computation of homoclinic orbits in delay and ordinary differential equations \cite{Samaey2002,Krauskopf2008}; see Ref.~\onlinecite{Homburg2010} for an extensive review of homoclinic orbits in the context of ordinary differential equations. 

In our case $\tau\in\mathbb{R}$, such that it is convenient to consider homoclinic orbits in the more general framework of delay differential equations of mixed type. Such equations naturally arise in the study of lattice dynamical systems, and a homoclinic solution to such an equation corresponds to a traveling pulse on the lattice\cite{Hupkes2009}.
A precise formulation of the problem relies on the existence of projection operators onto generalized stable and unstable eigenspaces of the corresponding steady state solution, together with the corresponding exponential dichotomies, which is guaranteed if the steady state is hyperbolic\cite{Mallet-Paret2001,Haerterich2002} (and the right handside of the linearized system does not vanish on an open interval of time). 

We here briefly review the steady states and their stability properties of Eqs.~(\ref{eq:G-def})--(\ref{eq:I-def}); see Ref.~\onlinecite{Terrien2017} for details. For certain parameter values there exist up to three equilibrium solutions $(G,Q,I)\in\mathbb{R}^3$ denoted $\mathsf{o},\mathsf{p},\mathsf{q}$. The laser off-state $\mathsf{o}=( A,B,0 )$ exists for all parameter values, and
$$p=\left(\frac{A}{1+I_-},\frac{B}{1+aI_-},I_-\right),\quad q=\left(\frac{A}{1+I_+},\frac{B}{1+aI_+},I_+\right),$$ where
\begin{align}
I_\pm= & \frac{ -aA+B+(1+a)(1-\kappa) \pm \sqrt{\Delta}}{2a(1-\kappa)}, \nonumber
\end{align}
exist a long as
\begin{align}
\Delta=&  \left( aA-B-(1+a)(1-\kappa) \right)^2 +\nonumber \\ 
& 4a(1-\kappa)\left(A-B-(1-\kappa)\right)\geq0.\nonumber
\end{align}
At $\Delta=0$, $\kappa=\kappa_\mathsf{F}$ and \textsf{p} and \textsf{q} undergo a fold bifurcation. At $\kappa=\kappa_\mathsf{T}$, \textsf{p} coincides with \textsf{o} in a transcritical bifurcation. For $\kappa>\kappa_\mathsf{T}$, the $I$-component of \textsf{p} is negative and \textsf{p} is no longer physically relevant. 

Let $\kappa<\kappa_\mathsf{T}$ in what follows. Given a steady state $\mathsf{r}\in\{\mathsf{o},\mathsf{p},\mathsf{q}\}$, a homoclinic orbit $\mathsf{l}(t)$ of Eqs.~(\ref{eq:G-def})--(\ref{eq:I-def}) is a solution $(G,Q,I)$ such that $(G(t),Q(t),I(t))\in W^{u}(\mathsf{r})$ as $t \to -\infty$ and\linebreak $(G(t),Q(t),I(t))\in W^{s}(\mathsf{r})$ as $t \to +\infty$.  The orbit $\mathsf{l_p}$ corresponding to \textsf{L\textsubscript{p}}  is homoclinic to \textsf{p}, which is a saddle for the considered parameter values; see Ref.~\onlinecite{Terrien2017} for details. This orbit corresponds to the homoclinic orbit in the case without feedback\cite{Dubbeldam1999,Huber2005}, which exists for $\tau\neq0$, and can be extended by continuation. 
We remark that, at $\tau=0$, countably many eigenvalues of the spectrum move through $-\infty$ to the right half plane; \textsf{p} remains a saddle. The continuation of the homoclinic orbit is numerically challenging; we rely on a numerical approach and continue \textsf{L\textsubscript{o}}  via a periodic orbit with sufficiently large period. 

The orbit $\mathsf{l_o}$ corresponding to \textsf{L\textsubscript{o}} is localized at \textsf{o}. Interestingly, for the parameter values considered in Fig.~\ref{fig:3}, the steady state \textsf{o} is a saddle, if and only if $\tau<0$ and, hence, $\mathsf{l_o}$ can only exists when $\tau<0$. This is due to the fact that the spectrum of the steady state \textsf{o} has infinitely many eigenvalues with positive real part if $\tau<0$. An overview of the spectral properties of \textsf{o} is contained in Appendix~\ref{sec:spec-o}. It is beyond the scope of this paper to rigorously proof the existence of $\mathsf{l_o}$; by using reappearance it is clear however, that the existence of a family of periodic solutions that is temporally localized about \textsf{o} with period $T(\tau)=\tau+\delta(\tau),$ and $\delta(\tau)\to\delta_0$ with $\delta_0/\tau\to\delta_0$ as $\tau\to\infty$ implies the existence of homoclinic orbit to \textsf{o} at $\tau=-\delta_0<0$. 

Such a homoclinic orbit, giving rise to periodic solutions, which are temporally localized about an exponentially stable steady state of a delay differential equation for $\tau\geq0$, is also called a temporal dissipative soliton\cite{Yanchuk2019}.  

\subsection{Stability of pulse trains corresponding to temporal dissipative solitons}\label{sec:TDS}

It has been shown recently that, for large values of the delay, the Floquet spectrum of the periodic pulse train (solid blue branches in Fig.~\ref{fig:3}(a) and (b)) corresponding to a temporal dissipative soliton splits into two distinct parts: the interface spectrum and the pseudo-continous spectrum. Whereas the interface spectrum consists of isolated Floquet multiplier as $\tau\to\infty$, the the pseudo-continuous spectrum can be approximated by a continuous curve, the so-called  asymptotic continuous spectrum; see Ref.~\onlinecite{Yanchuk2019} for details. As $\tau$ increases, this curve becomes densely filled with Floquet multipliers. 

We follow Ref.~\onlinecite{Yanchuk2019} to show that the pseudo-continuous spectrum of a $k$-pulse train (dark blue, cyan) is well approximated by the curve $\mu_k:\mathbb{R}\to\mathbb{C},$ where 
\begin{equation}\label{eq:acs}
\mu_k^k(\omega)=\frac{\kappa e^{i\omega\delta_0}}{i\omega-A+B+1},
\end{equation}
with an error of order $1/\tau$; the details are in Appendix~\ref{sec:pulse-train-spec}. Figure~\ref{fig4} shows the Floquet spectrum of the $1$-, $2$- and $3$-pulse train which coexist for $A=6.5$, $B=5.8$, $a=1.8$, $\gamma=0.04$, $\tau=4000$, and $\kappa=0.1$. We observe that the $1$-pulse train is stable and, as the number of pulses per delay interval $k$ increases, the modulus of the largest nontrivial multiplier increases like
$$\max_{\omega}\mu_k(\omega)=\sqrt[k]{\frac{\kappa}{|A-B-1|}}.$$  In addition, there are $k$ mulipliers close to the imaginary axis. Theses multipliers correspond to the translational Floquet eigendirections along which a $k$-pulse train is weakly attracting. 
Conditions $\kappa<\kappa_\mathsf{T}$ and (\ref{eq:crit-delay}) (equivalently (\ref{eq:crit-k})) imply that the 
Floquet multipliers approximated by the asymptotic continuous spectrum (and therefore, the asymptotic continuous spectrum itself) stay bounded away from the unit circle by a constant of order $1/\tau$. As $\kappa\to\kappa_\mathsf{T}$, the pseudo-continuous spectrum approaches the imaginary axis, because $\max_{\omega}\mu_k(\omega)\to1$. For $\kappa>\kappa_\mathsf{T}$, large delay instabilities can lead to complicated dynamics \cite{Klinshov2015a, Klinshov2015}. We note that periodic solutions for negative delay are always unstable. They exist as $C^\infty$ objects in the corresponding advanced differential equations, but a generic small perturbation along the periodic solutions diverges at least exponentially.

\begin{figure}[]
	\centering{t}
	\includegraphics[width=.6\linewidth]{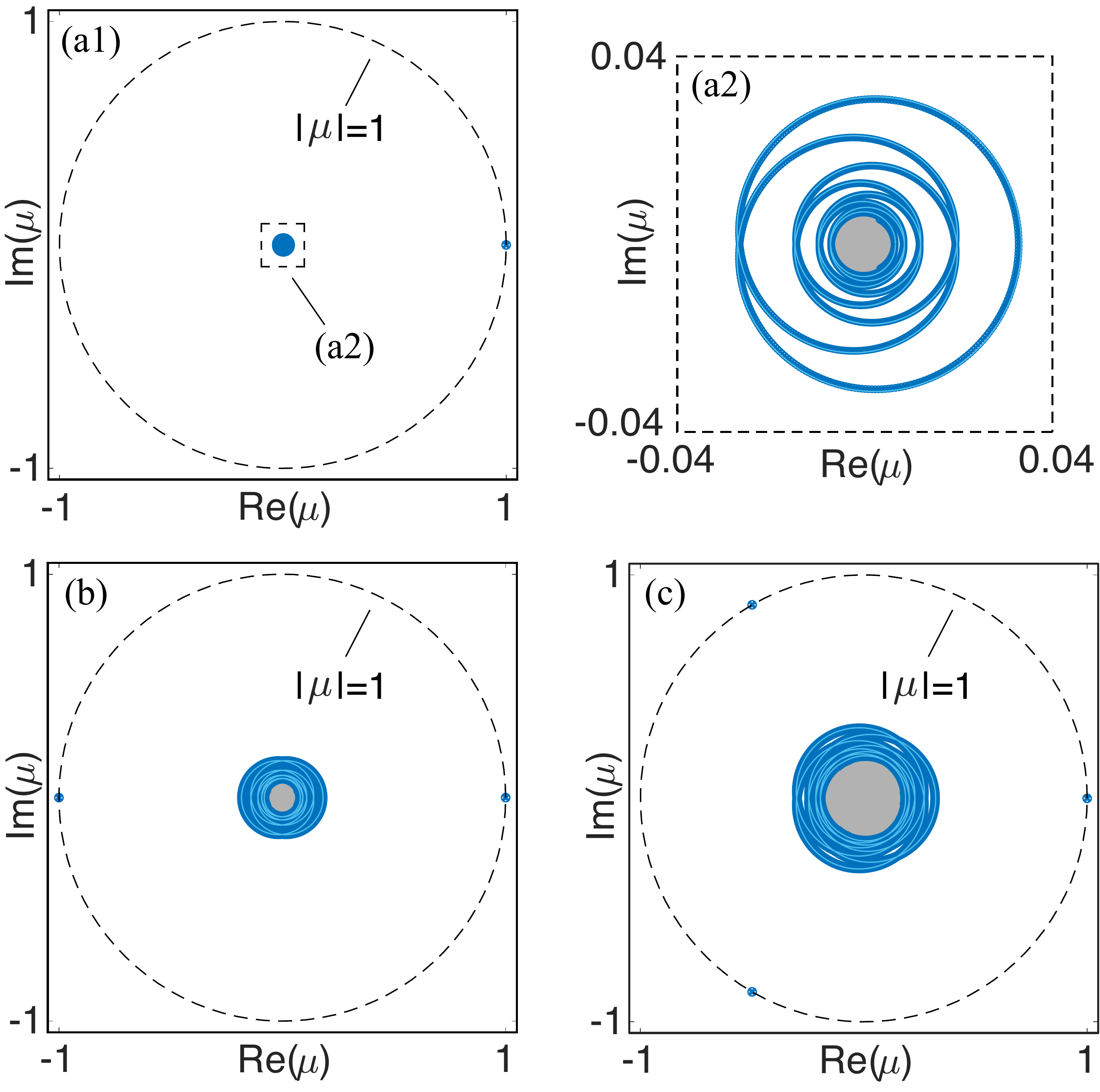}
	\caption{\label{fig4}
		Floquet spectrum of periodic pulse trains coexisting for parameter values $A=6.5$, $B=5.8$, $a=1.8$, $\gamma=0.04$, $\tau=4000$, and $\kappa=0.1$. Shown are the $2000$ Floquet multipliers (blue) largest in modulus and the asymptotic continuous spectrum (cyan). Floquet multipliers with smaller modulus have not been computed (gray regions). (a1) spectrum of the 1-pulse train; (a2) enlargement of (a1) close to $0$; (b) spectrum of the 2-pulse train; and (c) spectrum of the 3-pulse train.
	}
\end{figure}

\section[2-parameter continuation]{two-parameter bifurcation diagram}\label{sec:2parametercont}

Starting from the one-parameter bifurcation diagrams for different values of $\kappa$, we now investigate the qualitative changes in Eqs.~(\ref{eq:G-def})--(\ref{eq:I-def}) that lead to the disappearance of periodic pulse trains. We do so continuing the bifurcation points in Fig.~\ref{fig:3} as curves in the $(\tau,\kappa)$-parameter plane for different values of $A$. The resulting two-parameter bifurcation diagrams are shown in Figs.~\ref{fig:5}~and~\ref{fig:6}. 

\subsection*{Case A.$\quad$The laser without feedback is excitable}\label{sec:2parametercont-A}

\begin{figure*}[t]
	\centering{}\includegraphics[width=1.\linewidth]{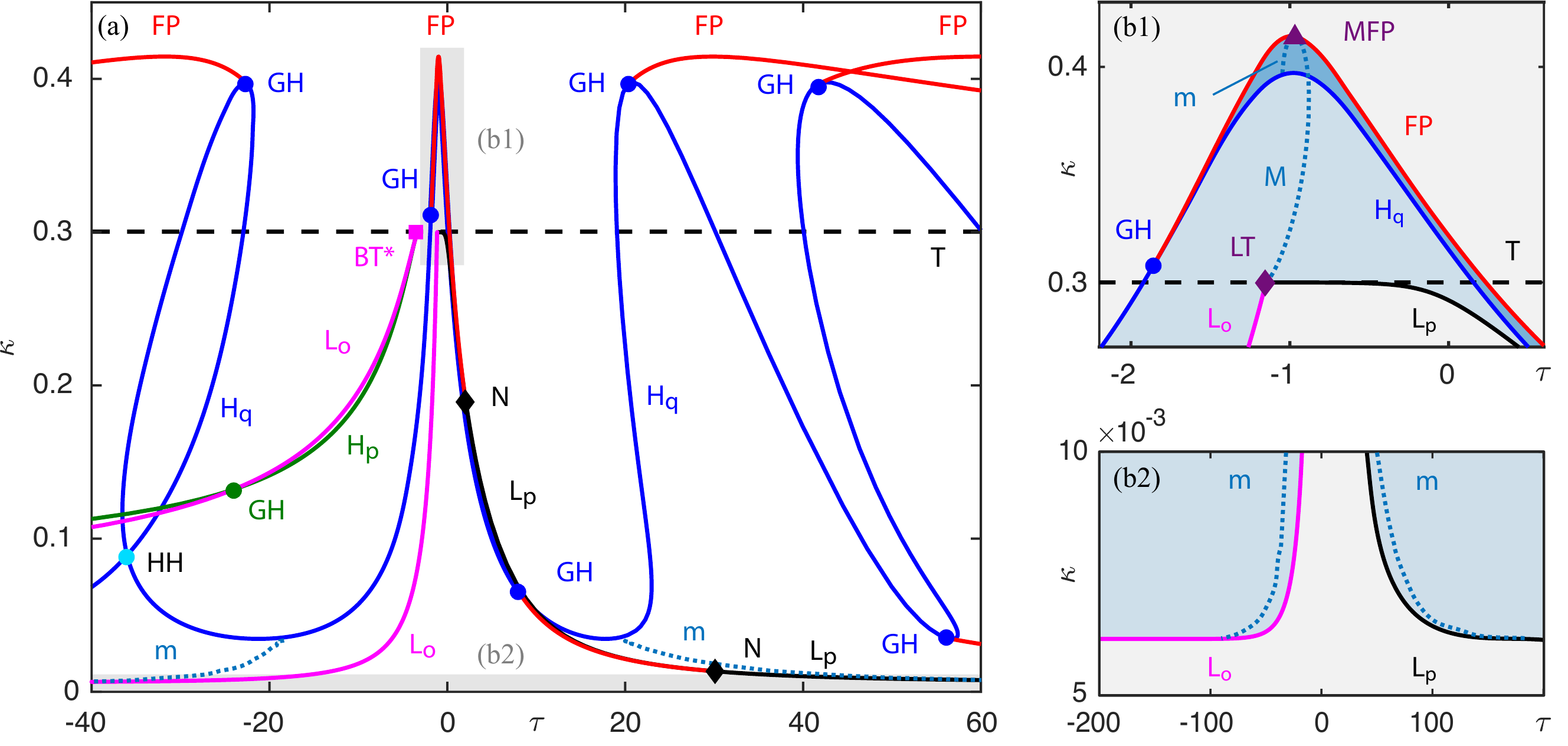}
	\caption{\label{fig:5}
		(Color online) Bifurcation diagram of Eqs.~(\ref{eq:G-def})--(\ref{eq:I-def}) in the $(\tau,\kappa)$-plane for $A=6.5$, $B=5.8$, $a=1.8$, and $\gamma_G=\gamma_Q=0.04$. (a) relevant part of the $(\tau,\kappa)$-plane; (b1)--(b2) enlargements of the respective gray shaded regions in (a) (and a larger range of $\tau$-values in (b2)). 
		Shown are curves of transcritical bifurcations \textsf{T} (dashed black), of Hopf bifurcations $\mathsf{H_p}$ (green) and $\mathsf{H_q}$ (blue), of homoclinic bifurcations \textsf{L\textsubscript{o}} (magenta) and \textsf{L\textsubscript{p}} (black), and of folds of periodic orbits \textsf{T} (red), as well as curves of periodic orbits with a local period minimum \textsf{m} or maximum \textsf{M} (dotted blue). 
		Additionally shown are points of generalized Hopf bifurcations \textsf{GH} (blue filled circle), of Neutral Saddle points N (filled diamond), a Hopf-Hopf bifurcation point \textsf{HH} (cyan filled circle), a Bogdanov-Takens bifurcation point \textsf{BT$^\ast$} (magenta filled square), a local maximum along a curve of fold periodic orbits \textsf{MFP} (purple triangle), as well as a homoclinic bifurcation point at nonhyperbolic equilibrium \textsf{LT}  (purple diamond).
		The light blue shaded regions in (b1)--(b2) indicate the existence of a periodic orbit; in the darker blue shaded region in (b1) periodic orbits coexist.
	}
\end{figure*}

First, we consider the parameter values $A=6.5$, $B=5.8$, $a=1.8$, $\gamma_G=\gamma_Q=0.04$ as in Fig.~\ref{fig:3}. This set up corresponds to the case when the laser without feedback is excitable. This exact set of parameters has been studied earlier by some of the authors\cite{Walker2008,Terrien2017}. We recover their results for $\tau\geq0$ and extend the bifurcation diagram to $\tau<0$, thereby allowing us to study the location of the homoclinic orbits \textsf{L\textsubscript{o}}  and \textsf{L\textsubscript{p}}; see Fig.~\ref{fig:5}.  A detailed description of the dynamics in the regions that are enclosed by these bifurcation curves for $\tau\geq 0$ can be found in Refs.~\onlinecite{Walker2008,Terrien2017}.

The bifurcation diagram in Fig.~\ref{fig:5}(a) is organized by Hopf bifurcations curves. For ease of reference, we denote Hopf bifurcation curves of the steady state \textsf{q} by \textsf{H\textsubscript{q}}  (blue), and of the steady state \textsf{p} by \textsf{H\textsubscript{p}}  (green), respectively. 
Along the curves \textsf{H\textsubscript{p}} and \textsf{H\textsubscript{q}}, we encounter generalized Hopf bifurcation points \textsf{GH} (filled circle). At a \textsf{GH} point, the Hopf bifurcation changes criticality (the direction in which a small perturbation of the parameter leads to the creation of a periodic orbit), which is generically accompanied by the emergence of a fold curve \textsf{FP} (red) of periodic orbits from the \textsf{GH} point.

We remark that, along the curve \textsf{H\textsubscript{p}}, the period of the bifurcating periodic orbit is large. The expected fold periodic orbit emanating from the point \textsf{GH} (green filled circle), therefore, is likely to have a large period as well. We were unable to find a fold curve of periodic orbits with period less than $T=10000$ near the point \textsf{GH} (green filled circle) in Fig.~\ref{fig:5}(a).

Two of the fold curves \textsf{FP} in Fig.~\ref{fig:5} terminate, respectively, in a neutral saddle bifurcation point \textsf{N} (filled black diamond), where the curve \textsf{FP} ends in the curve  \textsf{L\textsubscript{p}} (black) of homoclinic bifurcations. 
At $\kappa=0.3$, the curve \textsf{L\textsubscript{p}} coincides with the transcritical bifurcation curve \textsf{T} (black dashed) of steady states \textsf{o} and \textsf{p}. The point \textsf{LT} (filled purple diamond) is a homoclinic bifurcation point with a non-hyperbolic steady state, which is a codimension-two phenomenon. At the point \textsf{LT}, there emerges a curve of homoclinic orbits \textsf{L\textsubscript{o}} (magenta).

At $(\tau,\kappa)=(-3,0.3)$, we encounter a double zero eigenvalue along the curve \textsf{T}, which we refer to as a Bogdanov-Takens bifurcation point \textsf{BT$^\ast$} (magenta square, filled). This bifurcation point is different from a standard Bogdanov-Takens bifurcation in that it occurs on the invariant plane $I=0$; see Appendix~\ref{sec:spec-o} for details. The point \textsf{BT$^\ast$}  gives rise to the Hopf curve \textsf{H\textsubscript{q}} (green) together with an additional curve \textsf{L\textsubscript{o}} (magenta) of homoclinic orbits with small amplitude, which are not of particular interest here. Similarly, we observe that the Hopf curve \textsf{H\textsubscript{q}} self-intersects in a Hopf-Hopf bifurcation point \textsf{HH} (cyan filled circle), which is not studied further here.
 
The curves \textsf{m} and \textsf{M} (both light blue, dotted) correspond to periodic orbits having a local minimum and local maximum of the period with respect to $\tau$, respectively. As a reminder, we are interested in the location of the periodic orbits along \textsf{m} and \textsf{M} with respect to the homoclinic orbits \textsf{L\textsubscript{o}} and \textsf{L\textsubscript{p}}, as well as the position of the fold periodic orbits along \textsf{FP} shown in Fig~\ref{fig:3}. 

Starting from the horizontal slice $\kappa=0.01$ in Fig.~\ref{fig:5}(a) (see Fig.~\ref{fig:2}(a) for the one-parameter bifurcation diagram) and decreasing $\kappa$, we observe that the curves \textsf{L\textsubscript{o}}  and \textsf{L\textsubscript{p}}  level off at a certain value of $\kappa$; this corresponds to the minimum coupling strength $\kappa_{\min}\approx 0.006$ that allows for periodic pulse trains and, therefore, to the lower boundary of the self-excited oscillatory regime in Fig.~\ref{fig:2}; see Figure~\ref{fig:5}(b2) for an enlargement of (a) and a larger range of values $\tau$. In particular, the periodic orbits along \textsf{m} approach in the $(\tau,\kappa)$-plane a homoclinic orbit in \textsf{L\textsubscript{o}} or \textsf{L\textsubscript{p}}; thus, $T_{\min}$ grows beyond bound as $\kappa\searrow\kappa_{\min}$. The period along the curve \textsf{m} increases fast, and it is challenging to follow it for $T>1000$.  
Recall from Sec.~\ref{sec:reapperance} that, for fixed $\kappa$, fold periodic orbits accumulate at the respective periodic orbits with minimal or maximal period in parameter space as $\tau\to\infty$. Therefore, as $\kappa\searrow\kappa_{\min}$, there are countably many fold periodic orbits that undergo a bifurcation with period infinity at $\kappa=\kappa_{\min}$; for $\kappa<\kappa_{\min}$, there are no periodic pulse trains. 

We note that the periodic orbits along the curves \textsf{m} with the same period reappear exactly and, therefore, their profiles are identical. Let us assume that \textsf{L\textsubscript{o}}  (and \textsf{L\textsubscript{p}}) exists for all $\tau<200$ (and $\tau>200$) and the family of periodic orbits along \textsf{m} converges to \textsf{L\textsubscript{o}} (and, therefore, to \textsf{L\textsubscript{p}}  as $|\tau|\to\infty$). This implies a heteroclinic connection between \textsf{o} and \textsf{p}, as \textsf{L\textsubscript{o}}  and \textsf{L\textsubscript{p}}  coincide at $(\tau,\kappa)=(\pm\infty,\kappa_{\min})$.  Indeed, the profile of the homoclinic orbits along \textsf{L\textsubscript{o}}  and \textsf{L\textsubscript{p}}  indicate such a connection for $|\tau|\gg1$. 

Conversely, increasing $\kappa$ from $\kappa=0.01$ in Fig.~\ref{fig:5}(a), we observe that an additional fold curve \textsf{FP} of period orbits emerges from the curve \textsf{L\textsubscript{p}} at a neutral saddle point \textsf{N}. This curve \textsf{FP}, coincides with a Hopf bifurcation curve \textsf{H\textsubscript{q}} at a generalized Hopf bifurcation point \textsf{GH}. The remaining fold curves disappear in a similar way at \textsf{GH} points (only one is shown in Fig.~\ref{fig:5}(a)). 
At $\kappa=0.1$ all fold points \textsf{FP} have disappeared and all branches of periodic orbits terminate either in a Hopf bifurcation \textsf{H\textsubscript{q}}, or in a homoclinic orbit \textsf{L\textsubscript{o}}, or \textsf{L\textsubscript{p}}; see Fig.~\ref{fig:2}(b) for the one-parameter bifurcation diagram. 
When increasing $\kappa$ further in Fig.~\ref{fig:5}(a), another curve \textsf{FP} of fold periodic orbits arises from a neutral saddle point \textsf{N}; see Fig.~\ref{fig:2}(c) for the one-parameter bifurcation diagram for $\kappa=0.38$. The curve \textsf{FP} bounds the region of periodic orbits for small values of $\tau$; see Fig.~\ref{fig:5}(b1) for an enlargement. 
From the point \textsf{LT} in Fig.~\ref{fig:5}(b1) arises a curve \textsf{M} of local period maxima, resulting in the existence of countably many folds points of periodic orbits in the the branches of periodic orbits that reappear from \textsf{M}. Moreover, as $\kappa\searrow\kappa_\mathsf{T}$ the period of the orbits along \textsf{M} grows beyond bound, and so does the period of the orbits along a given curve \textsf{F} for $\kappa>\kappa_\mathsf{T}$. As a result, all curves \textsf{FP}, with $|\tau|\gg1$ and $\kappa>\kappa_\mathsf{T}$, approach \textsf{T} asymptotically as $|\tau|\to\infty$. In this way, the curve \textsf{T} bounds the region of periodic pulse trains with arbitrarily large period and, therefore, the self-excited oscillatory regime in Fig.~\ref{fig:2}.
	
We continue with the discussion of Fig.~\ref{fig:5}(b1). 
The fold curve \textsf{FP} displays a local maximum \textsf{MFP}, where the curves \textsf{FP}, \textsf{M} and \textsf{m} intersect. For $\kappa>\kappa_{\mathsf{MFP}}$, there are no periodic orbits. For $\kappa_T<\kappa<\kappa_{\mathsf{MFP}}$ periodic orbits have bounded period, and the orbit with the maximum period for fixed $\kappa$ is indicated by \textsf{M}. As a result, the boundary of the oscillatory regime $\kappa_{\max}$ is given by the $\kappa$-value of the point \textsf{MFP}. 
For a fixed value of $\kappa$ above the local maximum of the Hopf curve $\mathsf{H_q}$ in Fig.~\ref{fig:5}(b1), but smaller than the point \textsf{MFP}, there exists an isola of periodic orbits parametrized by $\tau$, which is bounded by fold points \textsf{FP} of periodic orbits. These isolas reappear for different values of the delay in Fig.~\ref{fig:5}(a).  The width of each of the $\tau$-intervals parametrizing these isolas shrinks to zero as $\kappa\nearrow\kappa_{\max}$.

The discussion above shows how $\kappa_{\min}$ and $\kappa_{\max}$ can be obtained by studying the two-parameter bifurcation diagram in the $(\tau,\kappa)$-plane. Recomputing the diagram for a large number of different values of $A$, one obtains Fig.~\ref{fig:2} by plotting the values of $\kappa_{\min}$ and $\kappa_{\max}$ with respect to the pump parameter $A$. For different values of $A$, the bifurcation diagram in the $(\tau,\kappa)$-plane changes of course. However, it remains effectively the same when $A$ is increased from the value $A=6.5$ in Fig.~\ref{fig:5}. 
More specifically, the curve \textsf{T} attains smaller and smaller values of $\kappa$, until the point when the laser without feedback undergoes a transcritical bifurcation. At this bifurcation point, the laser becomes self-oscillating and the value of $\kappa_{\min}$ is zero; see Fig.~\ref{fig:2}. 
As the value of $A$ increases, the maximum value of $\kappa$ along the curve \textsf{H\textsubscript{q}} decreases and so does the distance between curves \textsf{H\textsubscript{q}} and \textsf{FP} in Fig.~\ref{fig:5}(b1). At $A\approx 7.5$, the constant lasing state \textsf{q} of the laser without feedback undergoes a stabilizing Hopf-bifurcation. At this bifurcation point, $\kappa_{\max}=0,$ so that we do not observe any periodic solutions for larger values of $A$; see Fig.~\ref{fig:2}.

\subsection*{Case B.$\quad$ The laser without feedback is off}\label{sec:2parametercont-B} 

Decreasing the value of $A$, on the other hand, the bifurcation diagram in the $(\tau,\kappa)$-plane changes a lot.
Nevertheless, It is still possible to consistently define $\kappa_{\min}$ and $\kappa_{\max}$, as the minimal value of $\kappa$ along the curve \textsf{L\textsubscript{o}} and, the maximal value of $\kappa$ along a curve \textsf{FP}, respectively.

\begin{figure*}[]
	\includegraphics[width=1.\linewidth]{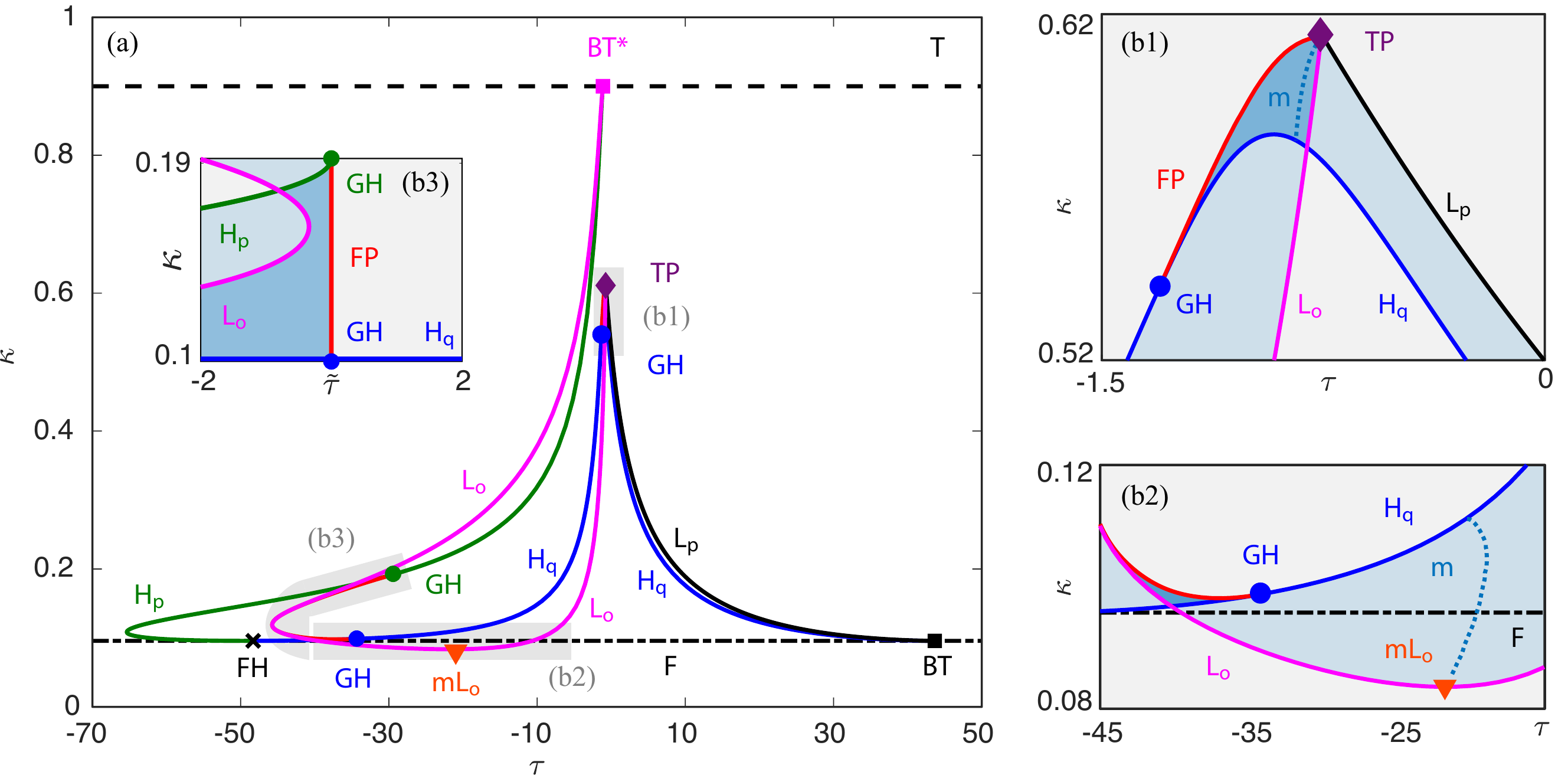}
	\caption{\label{fig:6}
		(Color online) Bifurcation diagram of Eqs.~(\ref{eq:G-def})--(\ref{eq:I-def}) in the $(\tau,\kappa)$-plane for $A=5.9$, $B=5.8$, $a=1.8$, and $\gamma_G=\gamma_Q=0.04$. (a) relevant part of the $(\tau,\kappa)$-plane; (b1)--(b3) enlargements of the respective gray shaded regions in (a). 
		Shown are curves of transcritical bifurcations \textsf{T} (dashed black), of fold (or saddle-node) bifurcations \textsf{F} (dashed dotted black), of Hopf bifurcations $\mathsf{H_p}$ (green) and $\mathsf{H_q}$ (blue), of homoclinic bifurcations \textsf{L\textsubscript{o}} (magenta) and \textsf{L\textsubscript{p}} (black), and of folds of periodic orbits \textsf{T} (red), as well as curves of periodic orbits with a local period minimum \textsf{m} or maximum \textsf{M} (dotted blue). 
		Additionally shown are points of generalized Hopf bifurcations \textsf{GH} (filled circle), of a fold-Hopf bifurcation point \textsf{FH} (black cross), of Bogdanov-Takens bifurcations \textsf{BT$^\ast$}  (black filled square) and \textsf{BT$^\ast$}  (magenta filled square), of a heteroclinic connection \textsf{TP} (purple diamond), as well as a local minimum  \textsf{mL\textsubscript{o}} (orange triangle) along the curve \textsf{L\textsubscript{o}}. 
		The light blue shaded regions in (b1)--(b3) indicate the existence of a periodic orbit; darker blue shaded regions in (b1)--(b3) indicate coexistence of periodic orbits. The horizontal axis in (b3) is shown relative a fold curve of periodic orbits between the two points of generalized Hopf bifurcations $\mathsf{GH}$ (blue) and $\mathsf{GH}$ (green) in the lower left of panel (a). 
	}
\end{figure*}

To illustrate this, we now show the bifurcation diagram in the $(\tau,\kappa)$-plane for $A=5.9$ (with all other parameter remaining the same); see Fig.~\ref{fig:6}. This configuration is interesting, because it corresponds to the case when the laser is off in the absence of feedback, but shows feedback-induced excitability for intermediate values of $\kappa$. The bifurcation diagram in Fig.~\ref{fig:6} is an extension of earlier results\cite{Terrien2017,Walker2008} to $\tau<0$, and we again focus on the location of the curves \textsf{m} and \textsf{M} relative to \textsf{L\textsubscript{o}}, \textsf{L\textsubscript{p}} and \textsf{FP}.

In Fig.~\ref{fig:6}(a), we now encounter a fold curve \textsf{F} (dashed dotted black) of steady states, along which \textsf{p} and \textsf{q} bifurcate; they exist when $\kappa>\kappa_F\approx0.1$.  Along the curve \textsf{F}, we find a fold-Hopf bifurcation point \textsf{FH} (dark cross), from which two Hopf curves \textsf{H\textsubscript{p}} (green) and \textsf{H\textsubscript{q}} (blue) emanate. The point \textsf{FH} may give rise to further curves of codimension-one bifurcations, which are not of immediate interest here. Following the curve \textsf{H\textsubscript{q}} to positive values of $\tau$, we find a (regular) Bogdanov-Takens bifurcation point \textsf{BT} (black filled square); see Fig.~\ref{fig:6}(b1) for an enlargement of (a) close to $\tau=0$. From the point \textsf{BT}, there emerges a curve \textsf{L\textsubscript{p}} of homoclinic orbits with small amplitude. Conversely, starting from the \textsf{FH} and increasing $\kappa$ along the curve \textsf{H\textsubscript{p}} , we find a Bogdanov-Takens bifurcation point \textsf{BT$^\ast$} (magenta filled square) on the transcritical bifurcation curve \textsf{T} (black dashed). A curve of homoclinic orbits \textsf{L\textsubscript{o}} with small amplitude emerges from \textsf{BT$^\ast$}. On the curves \textsf{H\textsubscript{p}} and \textsf{H\textsubscript{q}}, we again encounter generalized Hopf bifurcation points \textsf{GH} (filled circles) that give rise to fold curves \textsf{FP} (red) of periodic orbits.

We remark that the curves \textsf{L\textsubscript{o}}, \textsf{L\textsubscript{p}}, \textsf{H\textsubscript{p}}, and \textsf{H\textsubscript{q}} exist in a 
bounded interval of the delay; this is in contrast to the bifurcation diagram for $A=6.5$ in Fig.~\ref{fig:5}, where they are unbounded. In particular, the curve \textsf{L\textsubscript{o}} attains a minimum with respect to $\kappa$ at the point \textsf{mL\textsubscript{o}} (orange filled triangle), which corresponds to the minimum coupling strength $\kappa_{\min}$ of stable periodic pulse trains for $A=5.9$; see Fig.~\ref{fig:6}(b2) for an enlargement of panel (a). The period along the curve \textsf{m} of periodic solution with minimal period with respect to $\tau$ also grows beyond bound as $\kappa\searrow\kappa_{\min}$; for $\kappa<\kappa_{\min}$ there are no periodic pulse trains with $\tau\geq0$. Note that the critical value $\kappa_{\min}$ does no longer correspond to a bifurcation point.  

The curves \textsf{L\textsubscript{o}},  \textsf{L\textsubscript{p}} and a fold curve \textsf{F} coincide at the point \textsf{TP} (purple filled diamond), which corresponds to a connecting orbit between the steady states \textsf{o} and \textsf{p}; see Fig.~\ref{fig:6}(b1) for an enlargement of (a) close to $\tau=0$. This global bifurcation is similar to a T-point in ordinary differential equations\cite{Homburg2010,Knobloch2013}, which requires the existence of a generic connection between \textsf{o} and \textsf{p}.
Indeed, for the considered parameter values the dimension of the stable eigenspace of \textsf{o} is $\text{dim}W^{s}(\mathsf{o})=3$, and the codimension of the unstable eigenspace of \textsf{p} is $\text{codim}W^{u}(\mathsf{p})=2$, implying the possibility of a generic connection between \textsf{o} and \textsf{p}. The profile of the homoclinic orbits along the curves \textsf{L\textsubscript{o}} and \textsf{L\textsubscript{p}} further supports the existence of a connecting orbit at the point \textsf{TP}. For details on the spectrum of the steady state \textsf{o} see Appendix~\ref{sec:spec-o}. While we do not study this global bifurcation in detail here, it follows that for $\kappa>\kappa_\mathsf{TP}\approx0.62$ there are no stable periodic pulse trains with $\tau\geq0$.  Therefore, the boundary of the self-excited oscillatory regime $\kappa_{\max}$ for $A=5.9$ is given by the point \textsf{TP}. As $\kappa\nearrow\kappa_{\max}$, the minimum period $T_{\min}$ along the curve \textsf{m} in Fig.~\ref{fig:6}(b1) grows beyond bound.  

For ease of reference, we only show the relevant Hopf bifurcation curves \textsf{H\textsubscript{p}} and \textsf{H\textsubscript{q}} near $\tau=0$ in Fig.~\ref{fig:6}(a). Note however, that countably many curves \textsf{H\textsubscript{p}} and \textsf{H\textsubscript{q}} of Hopf bifurcations reappear under the mapping $\tau\mapsto\tau+2\pi/\omega,$ where $\omega$ is the imaginary part of critical eigenvalues of the respective equilibrium at the Hopf bifurcation points along \textsf{H\textsubscript{p}} and \textsf{H\textsubscript{q}}. Along theses curves, one encounters Hopf-Hopf bifurcation points and further points of generalized Hopf-bifurcations, each giving rise to additional curves of fold periodic orbits; see Ref.~\onlinecite{Walker2008}.
Recalling Sec.~\ref{sec:reapperance}, curves \textsf{m} reappear in a similar way, and the points along the curves of folds of periodic orbits above accumulate at points along these curves \textsf{m}. Therefore, as $\kappa\searrow\kappa_{\max}$, and respectively $\kappa\nearrow\kappa_{\min}$, the fold curves approach the line $\kappa=\kappa_{\max}$,  respectively $\kappa=\kappa_{\min}$, as $\tau\to\infty$.
The full bifurcation diagram for a finite range of positive delays $\tau$, can be found in Refs.~\onlinecite{Terrien2017,Walker2008}, together with a detailed analysis of the bifurcations of higher-codimension that are involved in the transition between the diagrams in Figs.~\ref{fig:5}(a)~and~\ref{fig:6}(a). 

As a side remark, along the Hopf bifurcation curve \textsf{H\textsubscript{p}} we encounter a generalized Hopf bifurcation point \textsf{GH}  (green filled dot), which gives rise to a fold curve \textsf{FP} of periodic orbits that ultimately disappears in another generalized Hopf bifurcation point \textsf{GH} (blue filled dot) on the curve \textsf{H\textsubscript{q}}. The curve \textsf{F} is very close in parameter space to the curve \textsf{L\textsubscript{o}}; Figure~\ref{fig:6}(b3) shows an enlargement of the corresponding parameter region; here the bifurcation curves are plotted in relative distance form the fold curve \textsf{FP} of periodic orbits. This shows that there are no further bifurcations in this small parameter region.

\section{Conclusions}\label{sec:discussion}


We studied the onset and termination of periodic pulse trains in the Yamada model for a semiconductor laser with saturable absorber and delayed optical feedback. A general bifurcation analysis of the model was given earlier by some of the authors\cite{Walker2008,Terrien2017}, and we formally extended this analysis to negative values of the feedback delay. 
With this approach, we were able to characterize the critical coupling strengths $\kappa_{\min}$ and $\kappa_{\max}$ that determine the onset and termination of periodic pulse trains, as well as the geometric mechanisms that are responsible for the corresponding dynamical changes. The main idea of the approach is to identify homoclinic and fold bifurcations of periodic orbits and to study the two-parameter bifurcation diagram in the feedback parameters, the delay $\tau$ and feedback strength $\kappa$.

More specifically, we numerically studied codimension-two points of homoclinic orbits and fold periodic orbits for different values of the pump parameter $A$, coupling strength $\kappa$, and feedback delay $\tau$. Projecting these points onto the $(A,\kappa)$-parameter plane revealed the curves $\kappa_{\min}$ and $\kappa_{\max}$ bounding the region of sustained self-pulsation. In this regime, the period of the periodic pulse train is largely determined by the length of the feedback loop, i.e. by the delay $\tau$. The regime of self-pulsation can be accessed (for sufficiently large coupling strength $\kappa_{\min}<\kappa<\kappa_{\max}$) even if the solitary laser is neither excitable nor oscillatory, well before the lasing threshold. The experimental validations of these predictions will be the subject of future research. 

As a secondary result, we employed the recently developed theory of temporal dissipative solitons\cite{Yanchuk2019} to infer upper estimates for the maximum number of equidistant pulses per delay interval that give rise to stable periodic pulse trains for fixed parameter values in the self-excited oscillatory region. Our results indicate that, for a given delay, there coexist a large number of stable periodic pulse trains that are weakly stable. Indeed, it has been observed in numerical simulations as well as in the experiment, that an arbitrary configuration of consecutive pulses convergences very slowly towards an equidistantly spaced pulse train\cite{Marino2017,Terrien2018a, Terrien2019}. The impact of noise on the stability of a given pulse train with $k$ pulses has not been considered here, but is worth investigating in the future. An earlier study showed, that, even in the oscillatory regime where there appear to be only a few coexisting stable pulse trains, the basins of attraction may be complicated Cantor-like sets \cite{Terrien2017}. Noise in the laser cavity with sufficiently large intensity may, therefore, be able to switch between pulse trains with different numbers of pulses. Similarly, one might expect such noise fluctuations to be able to trigger intermediate pulses, possibly causing the laser to accumulate the largest number of pulses possible for a given value of the feedback delay. 

For larger values of $\gamma_G$ and $\gamma_Q$, multiple modes of laser light start to interfere, thereby allowing for passive-mode locking inside the laser cavity\cite{Terrien2018, Vladimirov2005}. This alternative mechanism constitutes another way to generate pulsed laser light, which is not captured by the Yamada model.  Our method of identifying the critical coupling strength, however, applies in a similar way, and this will be a topic of future study. Interestingly, the corresponding delay differential equation model\cite{Vladimirov2005} allows for bound pulse trains\cite{Puzyrev2017} (a train of consecutive pulses followed by a longer pause), which correspond to $N$-homoclinic orbits. Such orbits and their bifurcations have been shown to organize the dynamical behavior of Lang-Kobayashi-type models of semiconductor lasers with optical injection\cite{Wieczorek2005}, and may facilitate our understanding of the complex dynamics in semiconductor lasers with optical feedback. Such $N$-homoclinic orbits have, however, not yet been observed in Eqs.~(\ref{eq:G-def})--(\ref{eq:I-def}). 

We want to emphasize that the presented methodology is not limited to the Yamda model with delay or to nonlinear optics, but is relevant for the rigorous analysis of delay-coupled excitable systems more generally. It may, furthermore, facilitate the understanding of delay-induced switched states; see Ref.~\onlinecite{Ruschel2019a} and references therein. On the other hand, there is a specific interest in the dynamics of coupled laser systems, because they are relevant for optical data storage and neuro-morphic information processing\cite{Bueno2018,Selmi2015,Prucnal2016}. 

\section*{Acknowledgments}
The authors thank Andrus Giraldo, Jan Sieber, Soizic Terrien, Matthias Wolfrum and Serhiy Yanchuk for valuable discussions related to this work.

\section*{Data Availability}
The data that support the findings of this study are available from the corresponding author
upon reasonable request.


\begin{thebibliography}{10}

\bibitem{Agrawal2002}
G.~P. Agrawal.
\newblock {\em {Fiber-Optic Communication Systems}}.
\newblock Wiley {\&} Sons, 2002.

\bibitem{Appeltant2011}
L.~Appeltant, M.~C. Soriano, G.~{Van Der Sande}, J.~Danckaert, S.~Massar,
  J.~Dambre, B.~Schrauwen, C.~R. Mirasso, and I.~Fischer.
\newblock {Information processing using a single dynamical node as complex
  system}.
\newblock {\em Nat. Commun.}, (2):468 |, 2011.

\bibitem{Argyris2005}
A.~Argyris, D.~Syvridis, L.~Larger, V.~Annovazzi-Lodi, P.~Colet, I.~Fischer,
  J.~Garc{\'{i}}a-Ojalvo, C.~R. Mirasso, L.~Pesquera, and K.~A. Shore.
\newblock {Chaos-based communications at high bit rates using commercial
  fibre-optic links}.
\newblock {\em Nature}, 438(7066):343--346, 2005.

\bibitem{Barbay2011}
S.~Barbay, R.~Kuszelewicz, and A.~M. Yacomotti.
\newblock {Excitability in a semiconductor laser with saturable absorber}.
\newblock {\em Opt. Lett.}, 36(23):4476, 2011.

\bibitem{Bueno2018}
J.~Bueno, S.~Maktoobi, L.~Froehly, I.~Fischer, M.~Jacquot, L.~Larger, and
  D.~Brunner.
\newblock {Reinforcement learning in a large-scale photonic recurrent neural
  network}.
\newblock {\em Optica}, 5(6):756, 2018.

\bibitem{Diekmann1995}
O.~Diekmann, S.~van Gils, S.~M. {Verduyn Lunel}, and H.-O. Walther.
\newblock {\em {Delay Equations, Functional-, Complex-, and Nonlinear
  Analysis}}.
\newblock Springer-Verlag, New York, 1995.

\bibitem{Dubbeldam1999}
J.~L. Dubbeldam and B.~Krauskopf.
\newblock {Self-pulsations of lasers with saturable absorber: Dynamics and
  bifurcations}.
\newblock {\em Opt. Commun.}, 159(4-6):325--338, 1999.

\bibitem{Dubbeldam1999a}
J.~L. Dubbeldam, B.~Krauskopf, and D.~Lenstra.
\newblock {Excitability and coherence resonance in lasers with saturable
  absorber}.
\newblock {\em Phys. Rev. E}, 60(6):6580--6588, 1999.

\bibitem{Engelborghs2002}
K.~Engelborghs, T.~Luzyanina, and D.~Roose.
\newblock {Numerical bifurcation analysis of delay differential equations using
  DDE-BIFTOOL}.
\newblock {\em ACM Trans. Math. Softw.}, 28(1):1--21, 2002.

\bibitem{Erneux1988}
T.~Erneux.
\newblock {Q-switching bifurcation in a laser with a saturable absorber}.
\newblock {\em J. Opt. Soc. Am. B}, 5(5):1063, 1988.

\bibitem{Garbin2017}
B.~Garbin, A.~Dolcemascolo, F.~Prati, J.~Javaloyes, G.~Tissoni, and S.~Barland.
\newblock {Refractory period of an excitable semiconductor laser with optical
  injection}.
\newblock {\em Phys. Rev. E}, 95:012214, 2017.

\bibitem{Gibson1974}
A.~F. Gibson.
\newblock {Laser-driven fusion}.
\newblock {\em Phys. Educ.}, 15(1):4--9, 1980.

\bibitem{Guo2013}
S.~Guo and J.~Wu.
\newblock {\em {Bifurcation Theory of Functional Differential Equations}},
  volume 184 of {\em Applied Mathematical Sciences}.
\newblock Springer New York, New York, NY, 2013.

\bibitem{Hale1993}
J.~K. Hale and S.~M.~V. Lunel.
\newblock {\em {Introduction to functional differential equations}}.
\newblock Springer, 1993.

\bibitem{Haerterich2002}
J.~H{\"{a}}rterich, B.~Sandstede, and A.~Scheel.
\newblock {Exponential dichotomies for linear non-autonomous functional
  differential equations of mixed type}.
\newblock {\em Indiana Univ. Math. J.}, 51(5):1081--1109, 2002.

\bibitem{Homburg2010}
A.~J. Homburg and B.~Sandstede.
\newblock {Homoclinic and heteroclinic bifurcations in vector fields}.
\newblock In B.~Hasselblatt, H.~W. Broer, and F.~Takens, editors, {\em Handb.
  Dyn. Syst.}, volume~3, pages 379--524. North-Holland, Amsterdam, 2010.

\bibitem{Huber2005}
A.~Huber and P.~Szmolyan.
\newblock {Geometric singular perturbation analysis of the Yamada model}.
\newblock {\em SIAM J. Appl. Dyn. Syst.}, 4(3):607--648, 2005.

\bibitem{Hupkes2009}
H.~J. Hupkes and S.~M. Lunel.
\newblock {Lin's method and homoclinic bifurcations for functional differential
  equations of mixed type}.
\newblock {\em Indiana Univ. Math. J.}, 58(6):2433--2487, 2009.

\bibitem{Izhikevich2005}
E.~M. Izhikevich.
\newblock {\em {Dynamical Systems in Neuroscience: The Geometry of Excitability
  and Bursting}}.
\newblock The MIT Press, 2007.

\bibitem{Jaurigue2015}
L.~Jaurigue, A.~Pimenov, D.~Rachinskii, E.~Sch{\"{o}}ll, K.~L{\"{u}}dge, and
  A.~G. Vladimirov.
\newblock {Timing jitter of passively-mode-locked semiconductor lasers subject
  to optical feedback: A semi-analytic approach}.
\newblock {\em Phys. Rev. A}, 92(5):053807 (2015), 2015.

\bibitem{Klinshov2015}
V.~Klinshov, L.~L{\"{u}}cken, D.~Shchapin, V.~Nekorkin, and S.~Yanchuk.
\newblock {Emergence and combinatorial accumulation of jittering regimes in
  spiking oscillators with delayed feedback}.
\newblock {\em Phys. Rev. E}, 92(4):45--47, 2015.

\bibitem{Klinshov2015a}
V.~Klinshov, L.~L{\"{u}}cken, D.~Shchapin, V.~Nekorkin, and S.~Yanchuk.
\newblock {Multistable jittering in oscillators with pulsatile delayed
  feedback}.
\newblock {\em Phys. Rev. Lett.}, 114(17):178103, 2015.

\bibitem{Knobloch2013}
J.~Knobloch, J.~S. Lamb, and K.~N. Webster.
\newblock {Using Lin's method to solve Bykov's problems}.
\newblock {\em J. Differ. Equ.}, 257(8):2984--3047, 2013.

\bibitem{Krauskopf2000}
B.~Krauskopf and D.~Lenstra.
\newblock {\em {Fundamental Issues of Nonlinear Laser Dynamics: Concepts,
  Mathematics, Physics, and Applications}}, volume 548.
\newblock American Institute of Physics, 2000.

\bibitem{Krauskopf2008}
B.~Krauskopf and T.~Rie{\ss}.
\newblock {A Lin's method approach to finding and continuing heteroclinic
  connections involving periodic orbits}.
\newblock {\em Nonlinearity}, 21(8):1655--1690, 2008.

\bibitem{Walker2008}
B.~Krauskopf and J.~J. Walker.
\newblock {Bifurcation Study of a Semiconductor Laser with Saturable Absorber
  and Delayed Optical Feedback}.
\newblock In {\em Nonlinear Laser Dyn. From Quantum Dots to Cryptogr.}, pages
  161--181. Wiley-VCH Verlag GmbH {\&} Co. KGaA, Weinheim, Germany, 2012.

\bibitem{Lin1986}
X.~B. Lin.
\newblock {Exponential dichotomies and homoclinic orbits in functional
  differential equations}.
\newblock {\em J. Differ. Equ.}, 63(2):227--254, 1986.

\bibitem{Lin1990}
X.~B. Lin.
\newblock {Using Melnikov's method to solve Silnikov's problems*}.
\newblock {\em Proc. R. Soc. Edinburgh Sect. A Math.}, 116(3-4):295--325, 1990.

\bibitem{Ludge2011}
K.~L{\"{u}}dge, editor.
\newblock {\em {Nonlinear Laser Dynamics: From Quantum Dots to Cryptography}}.
\newblock Wiley-VCH, 2012.

\bibitem{Mallet-Paret2001}
J.~Mallet-Paret and S.~V. Lunel.
\newblock {Exponential Dichotomies and Wiener-Hopf Factorizations for
  Mixed-Type Functional Differential Equations}.
\newblock {\em J. Differ. Equ.}, to appear, 2001.

\bibitem{Marino2017}
F.~Marino and G.~Giacomelli.
\newblock {Pseudo-spatial coherence resonance in an excitable laser with long
  delayed feedback}.
\newblock {\em Chaos}, 27(11), 2017.

\bibitem{Otto2012}
C.~Otto, K.~L{\"{u}}dge, A.~G. Vladimirov, M.~Wolfrum, and E.~Sch{\"{o}}ll.
\newblock {Delay-induced dynamics and jitter reduction of passively mode-locked
  semiconductor lasers subject to optical feedback}.
\newblock {\em New J. Phys.}, 14, 2012.

\bibitem{Otupiri2019}
R.~Otupiri, B.~Krauskopf, and N.~G.~R. Broderick.
\newblock {The Yamada model for a self-pulsing laser: bifurcation structure for
  non-identical decay times of gain and absorber}.
\newblock {\em Arxiv Prepr. 1911.01835}, 2019.

\bibitem{Pammi2019}
V.~A. Pammi, K.~Alfaro-Bittner, M.~G. Clerc, and S.~Barbay.
\newblock {Photonic Computing With Single and Coupled Spiking Micropillar
  Lasers}.
\newblock {\em IEEE J. Sel. Top. Quantum Electron.}, 26(1):1--7, 2019.

\bibitem{Prucnal2016}
P.~R. Prucnal, B.~J. Shastri, T.~{Ferreira de Lima}, M.~A. Nahmias, and A.~N.
  Tait.
\newblock {Recent progress in semiconductor excitable lasers for photonic spike
  processing}.
\newblock {\em Adv. Opt. Photonics}, 8(2):228, 2016.

\bibitem{Puzyrev2017}
D.~Puzyrev, A.~G. Vladimirov, A.~Pimenov, S.~V. Gurevich, and S.~Yanchuk.
\newblock {Bound Pulse Trains in Arrays of Coupled Spatially Extended Dynamical
  Systems}.
\newblock {\em Phys. Rev. Lett.}, 119(16):1--6, 2017.

\bibitem{RuschelYanchuk2017}
S.~Ruschel and S.~Yanchuk.
\newblock {Chaotic bursting in semiconductor lasers}.
\newblock {\em Chaos}, 27(11):114313, 2017.

\bibitem{Ruschel2019a}
S.~Ruschel and S.~Yanchuk.
\newblock {Delay-induced switched states in a slow-fast system}.
\newblock {\em Philos. Trans. R. Soc. A Math. Phys. Eng. Sci.}, 377(2153),
  2019.

\bibitem{Samaey2002}
G.~Samaey, K.~Engelborghs, and D.~Roose.
\newblock {Numerical computation of connecting orbits in delay differential
  equations}.
\newblock {\em Numer. Algorithms}, 30(3-4):335--352, 2002.

\bibitem{Selmi2015}
F.~Selmi, R.~Braive, G.~Beaudoin, I.~Sagnes, R.~Kuszelewicz, and S.~Barbay.
\newblock {Temporal summation in a neuromimetic micropillar laser}.
\newblock {\em Opt. Lett.}, 40(23):5690, 2015.

\bibitem{Sieber2014}
J.~Sieber, K.~Engelborghs, T.~Luzyanina, G.~Samaey, and D.~Roose.
\newblock {DDE-BIFTOOL Manual - Bifurcation analysis of delay differential
  equations}.
\newblock {\em arXiv}, 1406.7144, 2014.

\bibitem{Sieber2013a}
J.~Sieber, M.~Wolfrum, M.~Lichtner, and S.~Yanchuk.
\newblock {On the stability of periodic orbits in delay equations with large
  delay}.
\newblock {\em Discret. Contin. Dyn. Syst. Ser. A}, 33(7):3109--3134, 2013.

\bibitem{Soriano2013}
M.~C. Soriano, J.~Garc{\'{i}}a-Ojalvo, C.~R. Mirasso, and I.~Fischer.
\newblock {Complex photonics: Dynamics and applications of delay-coupled
  semiconductors lasers}.
\newblock {\em Rev. Mod. Phys.}, 85(1):421--470, 2013.

\bibitem{Steen2010}
W.~M. Steen and J.~Mazumder.
\newblock {\em {Laser material processing}}.
\newblock Springer, 4 edition, 2010.

\bibitem{Terrien2017}
S.~Terrien, B.~Krauskopf, and N.~G. Broderick.
\newblock {Bifurcation analysis of the yamada model for a pulsing semiconductor
  laser with saturable absorber and delayed optical feedback}.
\newblock {\em SIAM J. Appl. Dyn. Syst.}, 16(2):771--801, 2017.

\bibitem{Terrien2017a}
S.~Terrien, B.~Krauskopf, N.~G. Broderick, L.~Andr{\'{e}}oli, F.~Selmi,
  R.~Braive, G.~Beaudoin, I.~Sagnes, and S.~Barbay.
\newblock {Asymmetric noise sensitivity of pulse trains in an excitable
  microlaser with delayed optical feedback}.
\newblock {\em Phys. Rev. A}, 96(4):1--5, 2017.

\bibitem{Terrien2018}
S.~Terrien, B.~Krauskopf, N.~G. Broderick, L.~Jaurigue, and K.~L{\"{u}}dge.
\newblock {Q -switched pulsing lasers subject to delayed feedback: A model
  comparison}.
\newblock {\em Phys. Rev. A}, 98(4):1--10, 2018.

\bibitem{Terrien2018a}
S.~Terrien, B.~Krauskopf, N.~G.~R. Broderick, R.~Braive, G.~Beaudoin,
  I.~Sagnes, and S.~Barbay.
\newblock {Pulse train interaction and control in a microcavity laser with
  delayed optical feedback}.
\newblock {\em Opt. Lett.}, 43(13):3013, 2018.

\bibitem{Terrien2019}
S.~Terrien, V.~A. Pammi, N.~G.~R. Broderick, R.~Braive, G.~Beaudoin, I.~Sagnes,
  B.~Krauskopf, and S.~Barbay.
\newblock {Equalization of pulse timings in an excitable microlaser system with
  delay}.
\newblock {\em Arch. Prepr. 1907.11143}, 2019.

\bibitem{Ueno1985}
M.~Ueno and R.~Lang.
\newblock {Conditions for self-sustained pulsation and bistability in
  semiconductor lasers}.
\newblock {\em J. Appl. Phys.}, 58(4):1689--1692, 1985.

\bibitem{Vladimirov2005}
A.~G. Vladimirov and D.~Turaev.
\newblock {Model for passive mode locking in semiconductor lasers}.
\newblock {\em Phys. Rev. A}, 72(3):33808, 2005.

\bibitem{Wieczorek2005}
S.~Wieczorek, B.~Krauskopf, T.~B. Simpson, and D.~Lenstra.
\newblock {The dynamical complexity of optically injected semiconductor
  lasers}.
\newblock {\em Phys. Rep.}, 416:1--128, 2005.

\bibitem{Yamada1993}
M.~Yamada.
\newblock {A Theoretical Analysis of Self-Sustained Pulsation Phenomena in
  Narrow-Stripe Semiconductor Lasers}.
\newblock {\em IEEE J. Quantum Electron.}, 29(5):1330--1336, 1993.

\bibitem{YanchukGiacomelli2017}
S.~Yanchuk and G.~Giacomelli.
\newblock {Spatio-temporal phenomena in complex systems with time delays}.
\newblock {\em J. Phys. A Math. Theor.}, 50(10):103001, 2017.

\bibitem{Yanchuk2009}
S.~Yanchuk and P.~Perlikowski.
\newblock {Delay and periodicity}.
\newblock {\em Phys. Rev. E}, 79(4):046221, 2009.

\bibitem{Yanchuk2019}
S.~Yanchuk, S.~Ruschel, J.~Sieber, and M.~Wolfrum.
\newblock {Temporal Dissipative Solitons in Time-Delay Feedback Systems}.
\newblock {\em Phys. Rev. Lett.}, 123(5):053901, 2019.

\end{thebibliography}

\appendix
\clearpage

\section{Spectrum of the off-state \textsf{o}}\label{sec:spec-o}

This section contains a characterization of the spectrum of the steady state \textsf{o} of Eqs.~(\ref{eq:G-def})--(\ref{eq:I-def}). The characteristic equation at the steady state \textsf{o} reads 
\begin{equation}\label{eq:DDE-general-ce-0}
0=\left(\lambda+\gamma_G\right)\left(\lambda+ \gamma_Q\right)\left(-\lambda+ A-B-1+ \kappa e^{- \tau\lambda}\right). 
\end{equation}
A solution $\lambda\in\mathbb{C}$ to Eq.~(\ref{eq:DDE-general-ce-0}) is called an eigenvalue of the steady state \textsf{o}; the entirety of eigenvalues is called the spectrum. General information on the spectrum of a steady state of delay differential equations can be reviewed in classic textbooks\cite{Hale1993, Diekmann1995, Guo2013}. We note that for $\tau>0$, the real parts of eigenvalues in the spectrum accumulate at $-\infty$, whereas for $\tau<0$, the situation is reversed and the real parts of eigenvalues in the spectrum accumulate at $+\infty$. The following mental picture is very helpful. As $\tau$ passes through $0$ from positive to negative values, countably many eigenvalue with negative real part pass through $-\infty+i\mathbb{R}$ and appear from $+\infty+i\mathbb{R}$ with positive real part. 

We continue with the disucssion of Eq.~(\ref{eq:DDE-general-ce-0}). The two eigenvalues $\lambda_{1,2}=-\gamma_{G,Q}$ correspond to the attracting eigendirections $(1,0,0)^T$ and $(0,1,0)$ that span the invariant plane $I=0$. As a remark, a homoclinic orbit of \textsf{o} cannot arrive or leave (in backward time) via one of these directions, because $I=0$ is invariant. We can investigate the solutions of the remaining factor of Eq.~(\ref{eq:DDE-general-ce-0})  independently. These solutions correspond to the eigendirection $(0,0,1)$ perpendicular to the plane $I=0$. In particular, we can derive conditions, under which a pair of conjugated eigenvalues crosses the imaginary axis (which generically corresponds to a Hopf bifurcation on the nonlinear level). Under the assumption $\lambda =i\omega$, we can parametrically solve Eq.~(\ref{eq:DDE-general-ce-0}) for the curves
$$
\kappa(\omega) = \pm\sqrt{\omega^2+( A -  B -1)^2},
$$
$$
\tau(\omega) = -\frac{1}{\omega}\arg\left(\frac{i\omega - A + B +1}{\kappa(\omega)}\right)+\frac{2\pi k}{\omega}, 
$$
in the $(\tau,\kappa)$-plane, where $k\in\mathbb{Z}$ and we require that $\kappa>|A-B-1|$. If $\kappa<|A-B-1|$ no eigenvalues can cross the unit circle. Assume otherwise, then a pair of eigenvalues crosses the unit circle at $(\kappa,\tau)=(\kappa(\omega),\tau(\omega))$; moreover, for fixed $\kappa$, the number of pairs of eigenvalues that cross the imaginary axis grows beyond bound as $|\tau|\to\infty.$ 
 
\begin{figure}[]
	\centering
	\makegapedcells
	\centering
	\resizebox{.7\linewidth}{!}{
		\begin{tabular}{c|c|c|c|c}
			\multirow{2}{*}{~} & \multicolumn{2}{c|}{$A<B+1$} & %
			\multicolumn{2}{c}{$A>B+1$} \\
			& $\kappa<-(A-B-1)$ & $\kappa>-(A-B-1)$ & $\kappa>A-B-1$ & $\kappa<A-B-1$ \\
			\cline{1-5}
			\rotatebox{90}{$\quad~ \tau>0 $}\mbox{~}  & \mbox{~}\includegraphics[width=.2\columnwidth]{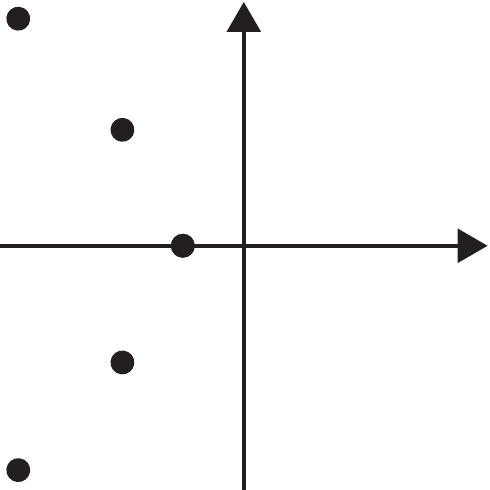} \mbox{\negthickspace}&
			\mbox{~}\includegraphics[width=.2\columnwidth]{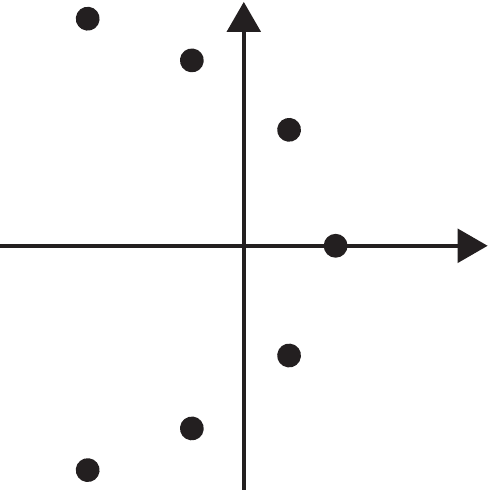}\mbox{\negthickspace}
			&\mbox{~}\includegraphics[width=.2\columnwidth]{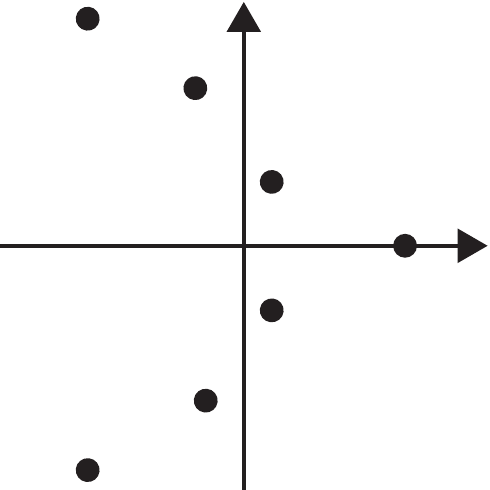}\mbox{\negthickspace} & \mbox{~}\includegraphics[width=.2\columnwidth]{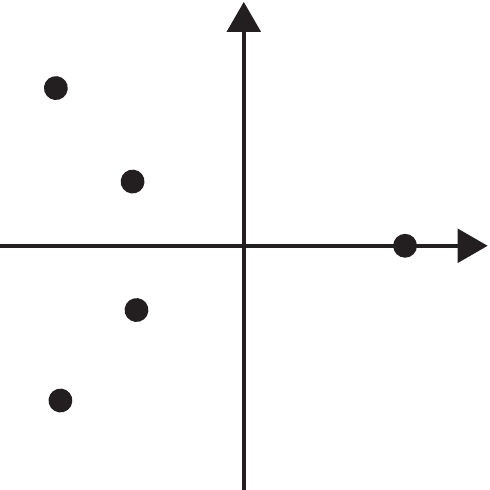}\mbox{\negthickspace} \\
			\cline{1-5}
			\rotatebox{90}{$\quad~ \tau<0 $}\mbox{~} & \mbox{~}\includegraphics[width=.2\columnwidth]{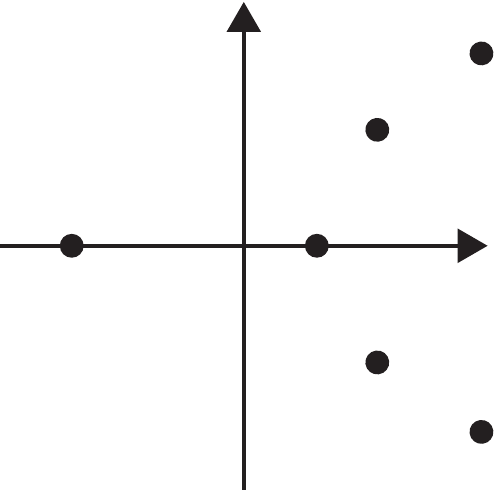}\mbox{~\negthickspace} & \mbox{~}\includegraphics[width=.2\columnwidth]{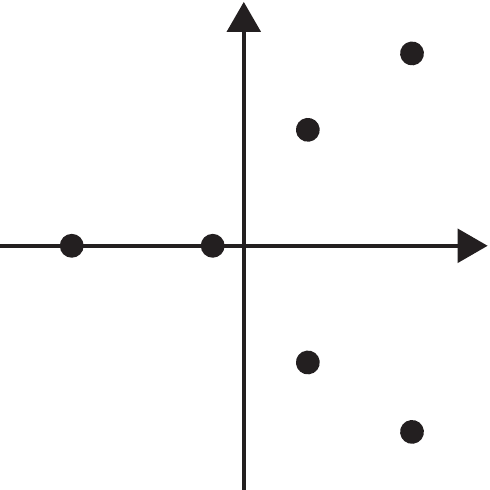}\mbox{\negthickspace}
			&\mbox{~}\includegraphics[width=.2\columnwidth]{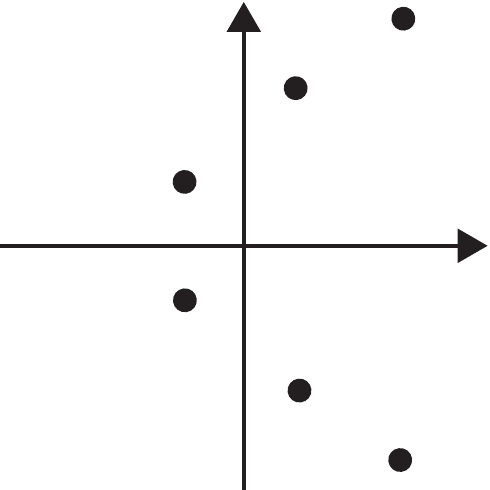}\mbox{\negthickspace} & \mbox{~}\includegraphics[width=.2\columnwidth]{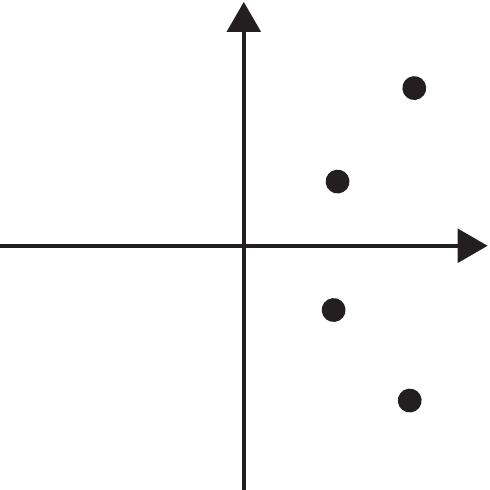}\mbox{\negthickspace} \\
		\end{tabular}
	}
	\caption{Sketch of the spectrum of $\mathsf{o}=(A,B,0)$ in the direction $(0,0,1)$ for sufficiently large $|\tau|$.  \label{tab:spec-0}}
\end{figure}

On the other hand, it follows from straightforward computation that there are no real solutions to Eq.~(\ref{eq:DDE-general-ce-0}), additional to $-\gamma_G$ and $-\gamma_Q$, if and only if $A>B+1$.  With this information at hand, one can completely characterize the spectrum of \textsf{o} in the direction perpendicular to $I=0$; this is illustrated in Fig.~\ref{tab:spec-0}. Throughout this work, we have considered parameter values satisfying $A<B+1$. In Figs.~\ref{fig:5} and \ref{fig:6}, we have shown that the homoclinic orbits along the curve \textsf{L\textsubscript{o}} exist for $\kappa>\kappa_\mathsf{T}=-(A-B-1)$ and $\tau<0$, corresponding to the case when the steady state \textsf{o} is a saddle with three-dimensional stable manifold (counting the stable directions on the invariant plane $I=0$). As $\kappa\to\kappa_\mathsf{T}$ and $\tau<0$, a real positive eigenvalue crosses the imaginary axis; this corresponds to the bifurcation point \textsf{LT} in Fig.~\ref{fig:5}(b1).
Conversely, for $\kappa\leq\kappa_\mathsf{T}$ and $\tau>0$ the steady state \textsf{o} is exponentially stable, and no homoclinic connection is possible. 

As a side remark, we find the location $$( \tau, \kappa)=\left(\frac{1}{A - B -1}, -(A - B -1)\right)$$ of the Takens-Bogdanov bifurcation point \textsf{BT$^\ast$} shown in Figs.~\ref{fig:5} and~\ref{fig:6} directly from Eq.~(\ref{eq:DDE-general-ce-0}).
  
\renewcommand{\arraystretch}{1.5}%
\renewcommand{\arraystretch}{1.}%

\section{Asymptotic continuous spectrum of stable periodic pulse trains}\label{sec:pulse-train-spec}
This sections contains the derivation of expression (\ref{eq:acs}) for the asymptotic continuous spectra shown in Fig.~\ref{fig4}. Consider a periodic solution $\tilde x(t)=(\tilde G(t),\tilde Q(t),\tilde I(t))$ of Eq.~(\ref{eq:G-def})--(\ref{eq:I-def}) with period $T$.  
The Floquet spectrum of $\tilde x$ is given by nontrivial (and normalized) solutions  $(\mu,y)\in\mathbb{C}\times C([-T,0],\mathbb{C}^n)$ to the equation 
\begin{equation}
y^\prime (t) = M_1(\tilde x(t))y(t) + M_2 y(t-\tau),\label{eq:Floquet}
\end{equation}
satisfying 
\begin{equation}
y(t)=\mu y(t+T),\label{eq:bddcon}
\end{equation}
for all $t$, where $M_1(\tilde x(t))$ and $M_2$ are given by the $T$-periodic coefficient matrices
$$
M_1(\tilde x(t)) = \begin{pmatrix}
-\gamma_Q (1+\tilde I(t)) & 0& -\gamma_Q \tilde G(t)\\
0 & -\gamma_Q (1+a \tilde I(t)) & -\gamma_Q a \tilde Q (t)\\
-\tilde I(t)  & -   \tilde I(t)  & \tilde G(t)- \tilde Q(t)-1
\end{pmatrix}
$$
and 
$$
M_2 = \begin{pmatrix}
0& 0& 0\\
0& 0& 0\\
0& 0& \kappa
\end{pmatrix}.
$$
The parameter $\mu$ is called the Floquet multiplier, and $y$ is called the eigensolution of $\mu$.  

More specifically, let $\tilde x$ be a $k$-pulse train in the self-excited oscillatory regime (see Fig.~\ref{fig:2}) with period $T=(\tau+\delta)/k$ (where $\delta>0$ is the drift of the pulse train, i.e. the response time of the laser to produce a subsequent pulse after delay time $\tau$) and $\delta=\delta(\tau)\to\delta_0$ as $\tau\to\infty$. As a reminder, we denote $(\tau,\kappa)=(-\delta_0,\kappa)$ the point at which there is a homoclinic bifurcation along the branch \textsf{L\textsubscript{o}} for fixed $\kappa$; see Figs.~\ref{fig:5} and~\ref{fig:6}. Henceforth, $\mu_k$ is a Floquet multiplier corresponding to $\tilde x$. 
Then, inserting Eq.~(\ref{eq:bddcon}) into Eq.~(\ref{eq:Floquet}) $k$-times, one obtains an equivalent formulation of the Floquet problem, where
\begin{equation}
y^\prime (t) = M_1(\tilde x(t)) + \mu_k^k M_2 y(t+\delta),\label{eq:Floquet-0}
\end{equation}
and $$y(t)=\mu_k y(t+T)$$ for all $t$.
As $\tau\to\infty$, the period of the pulse train grows beyond bound; equivalently, $\delta\to\delta_0$.

Now for sufficiently large values of $\tau$, the Floquet spectrum of $\tilde x$ splits into two parts: the interface spectrum and the pseudo-continuous spectrum \cite{Yanchuk2019}. For sufficiently large $\tau$, the Floquet multipliers in the pseudo-continuous spectrum accumulate along a continuous curve, the so called asymptotic continuous spectrum. Following Ref.~\onlinecite{Yanchuk2019}, the asymptotic continuous spectrum can be determined by evaluating Eq.~(\ref{eq:Floquet-0}) at $\tilde  x = \mathsf{o}$ for $\delta=\delta_0$, and choosing the spectral parameter $\mu_k$ such that the resulting steady state spectrum is non-hyperbolic, i.e.
\begin{equation}
0 = \det(-i\omega I + M_1(\mbox{\textsf{o}}) + \mu_k^k  M_2 e^{i\omega\delta_0})\label{eq:Floquet-acs}
\end{equation}
for some $\omega\in\mathbb{R}$. Solving (\ref{eq:Floquet-acs}) parametrically for $\mu_k^k=\mu_k^k(\omega)$ gives expression (\ref{eq:acs}).
\end{document}